\documentclass[12pt]{article}
\usepackage{amsfonts}
\usepackage{graphicx}
\usepackage[all,arc]{xy} \input xy

\newtheorem{thm}{Theorem}[section]
\newtheorem{lemma}[thm]{Lemma}
\newtheorem{cor}[thm]{Corollary}

\newenvironment{remark}{\par\medskip\noindent{\bf Remark.\ }}{\par\smallskip}
\newcommand{\proof
}{\par\medskip\noindent {\bf Proof.\ \ }}

\newcommand{\tah}[1]{\mbox{th}\left( #1 \right)}
\newcommand{\sh}[1]{\mbox{sh}\left( #1 \right)}
\newcommand{\ch}[1]{\mbox{ch}\left( #1 \right)}
\newcommand{\be}{\begin{equation}}
\newcommand{\ee}{\end{equation}}
\newcommand{\openbox}{\leavevmode
  \hbox to8pt{\hfil\vrule\vbox to6pt{\hrule width6pt\vfil\hrule}\vrule}}

\newcommand{\qed}{\hbox to5pt{ } \hfill \openbox\bigskip\medskip}

\newcommand{\Hi}{\mathbb H}

\newcommand{\R}{\mathbb R}

\title{The Minkowskian planar 4R mechanism}
\author{G{\'a}bor Heged{\"u}s \\
{\normalsize Johann Radon Institute for Computational and Applied Mathematics} \\
Brian Moore \\
\normalsize Advanced Telecommunications Research Institute International
}
\begin{document}

\maketitle

\begin{abstract}
We characterize and classify completely 
the planar 4R
closed chain working on the Minkowskian plane. Our work would open 
a new research direction 
in the theory of geometric designs:
the classification and characterization of the 
geometric design of linkages working in 
non--Euclidean spaces.
\end{abstract}
\section{Introduction}

First we recall here some preliminary definitions and results
from the geometry of linkages.

{\em A linkage} is a collection of interconnected components,
individually called {\em links}. The {\em joint} is the physical connection 
between two links. 

In our article we consider only the {\em revolute} joint, which can be viewed as 
constructed from the rotary hinge. We denote the revolute joint by R.

The revolute joint allows {\em one-degree-of-freedom} movement between 
the two links that it connects.
The configuration variable for a hinge is the angle 
measured around its axis 
between the two bodies.

Of course we can form linkages from other joints, for example 
the universal joint,
the ball-in-socket and the {\em prismatic} joint. 

%But there are other more exotic joints that can be used to construct
%linkages, such as the screw, or helix joint.

The {\em generic mobility} of the system is the number of independent 
parameters such as the joint angles 
that are needed to specify the configuration of the linkage.

On the other hand it can be shown that this is the dimension of the 
{\em configuration space} of the system.

\begin{figure}[htp]
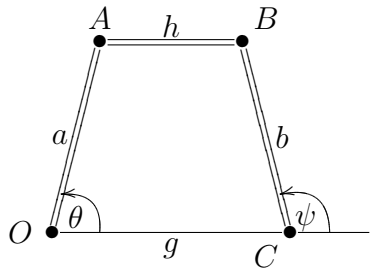

    \centering 
% sajat planar 4R
\[
\xy
0; <3pt,0pt>:
(0,0)*{\bullet}="O"; 
(30,0)*{\bullet}="C";
(6,24)*{\bullet}="A"; 
(24,24)*{\bullet}="B";
(40,0)*{}="Y";
(35,0)*{}="K";
(21,36)*{}="T";
(22.5,30)*{}="S";
"O"+(-4,0)*{O}; 
"C"+(-3,-3)*{C};
"A"+(0,3)*{A};
"B"+(3,3)*{B};
"O"; "A"; **\dir2{-}; ?(.8)*{}="X2"; ?(.2)*{}="X3";
"O"; "C"; **\dir{-}; ?(.8)*{}="X1";  ?(.2)*{}="X8";
"A"; "B"; **\dir{=}; ?(.8)*{}="X4"; ?(.3)*{}="X5"; 
"B"; "C"; **\dir{=};  ?(.8)*{}="X6"; ?(.2)*{}="X7";
"C"; "Y"; **\dir{-};
%"B"; "T"; **\dir{-};
%"X1"; "X2"; **\crv{(5,5)};
{\ar@/_0.5pc/ "X1";"X2"};
%"X3"; "X4"; **\crv{(11,19)};
%"X5"; "X6"; **\crv{(19,19)};
%"X7"; "X8"; **\crv{(25,5)};
%"X5"; "S"; **\crv{(17,28)};
%"X7"; "K"; **\crv{(34,5)};
{\ar@/_0.6pc/ "K";"X7"};
"O"+(3,2)*{\theta};
"C"+(2,2)*{\psi};
%"B"+(-3,3)*{\zeta};
%"A"+(2,-2)*{\phi};
(15,-2)*{g};
(1,12)*{a};
(15,26)*{h};
(29,12)*{b};
\endxy
\]
\caption[The planar 4R linkage]{The planar 4R linkage}
\label{planar_4R_linkage}
\end{figure}

A {\em planar} linkage has the property that all of its links 
move in parallel planes. 
%Most linkages which we find in practice 
%are planar linkages. 
We are interested here in the four--bar linkage,
which is a closed chain formed by four links and four joints. Figure
 \ref{planar_4R_linkage} is an example of a planar 4R closed chain.

The usual model of geometric designs are working in the
Euclidean space. 
In this article we would like to characterize and classify completely 
the planar 4R
closed chain working on the Minkowskian plane. This work would open 
a new research direction 
in the theory of geometric designs:
the classification and characterization of the 
geometric design of linkages working in 
non--Euclidean spaces.

In Chapter 2 we collected the preliminary definitions and results
 about the Minkowskian plane 
and hyperbolic trigonometry. In Chapter 3 we describe our main results: the position 
analysis and the classification of the Minkowskian 4R planar linkage.  We give also an 
exact formula for the coupler curve of this system and 
compute the transmission and 
coupler angles.

\section{Preliminaries}
\subsection{The Minkowskian plane and the double numbers}

First we give an algebraic description of the Minkowskian plane.

In analogy with the complex number system, the system of {\em double numbers}
\footnote{these numbers are called to split--complex numbers, too} 
can be introduced:

$$
\Hi:=\{x+jy:~ x,y\in {\R},\ j^2=1\}
$$  

Here $j$ is the {\em double imaginary unit} and $x$ and $y$ are respectively 
called the {\em real} and the {\em unipotent} parts of the double 
number $z=x+yj$.

It follows that multiplication in $\Hi$ is defined by 
$(x+yj)(r+sj)=(xr+ys)+j(xs+yr)$.

It is known that the complex numbers are related to 
the Euclidean geometry. Similarly
the double system of numbers serve as coordinates in
the Minkowskian plane (space-time geometry, see \cite{FG}).

The {\em hyperbolic conjugate} $\overline{z}$ of  $z=x+yj$ is defined by
$\overline{z}=x-yj$. 

The {\em hyperbolic scalar product} is given by
$$
\langle z,w\rangle:= Re(z\overline{w})=xu-yv,
$$
where $z=x+yj$ and $w=u+jv$. We say that the double numbers 
$z$ and $w$ are {\em double--orthogonal}, if 
$\langle z,w\rangle =0$.

We define the {\em hyperbolic modulus} of $z=x+yj$ by
\begin{equation}
\parallel z  {\parallel}_h:=
\sqrt{|\langle z,z\rangle |}
=\sqrt{|z \overline{z}|}=\sqrt{| x^2-y^2 |}\geq 0.
\end{equation}

We can consider this modulus  
as the {\em hyperbolic distance} of the point $z$ from the origin.  

It can be shown that this modulus is the Lorentz 
invariant of two dimensional special relativity, see \cite{Y}.

Note that the points $z\neq 0$ on the lines $y=x$ and $y=-x$ are
isotropic. This means that they are nonzero vectors with 
$\parallel z  \parallel_h =0$.

\begin{figure}[htp]
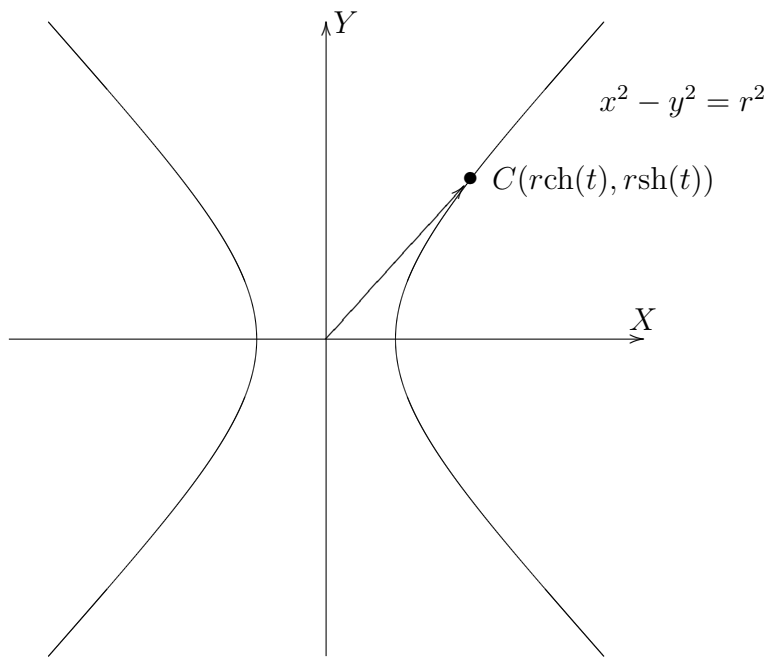

    \centering 
\[\xy
  0; <15pt,0pt>: % vector basis
 (-7,-8)*{}="-x"; (7,-8)*{}="x"; (-7,8)*{}="-y"; (7,8)*{}="y"; (0,0)*{}="O"; 
 (8,0)*{}="a";
(-8,0)*{}="-a";
(0,8)*{}="b";
(0,-8)*{}="-b";
"a"+(0,0.5)*{X};
"b"+(0.5,0)*{Y};
"-y"; "-x" **\crv{(0,0)&}  ?(.8)*{}="D";
 "x"; "y" **\crv{(0,0)&} ?(.8)*{\bullet}="C";
{\ar "O"; "C"}; 
{\ar "a"; "-a"; };
{\ar "b"; "-b"; };
(9,6)*{x^2-y^2=r^2};
(7,4)*{C(r \mbox{ch}(t),r \mbox{sh}(t))};
\endxy\]
\caption[hyperbola]{The Minkowskian circle of radius $r$ centered at the origin }
\label{hyperbola}
\end{figure}

For $r\in {\R}^+$ the {\em Minkowskian circle} of radius $r$ centered at the origin 
in $\Hi$ is defined by
$$
\{(x,y)\in {\R}:~ x^2-y^2=r^2 \}
$$ 
Clearly this set is the set of all points in the Minkowskian plane
 that satisfy 
the equation $\parallel z  \parallel_h^2 =r^2$ (see Figure \ref{hyperbola}).

We can remind here for the following analog of Euler's formula
for the double numbers $z=x+yj$:
$$
z=re^{j\phi}=r(\ch\phi+j\sh\phi)
$$
where $x^2-y^2>0$ and $y>0$.

\subsection{Hyperbolic trigonometry}
 
Let $L:=L^2$ denote the vector space ${\R}^2$ provided with 
the hyperbolic scalar product.

We denote by $SO^+(1,1)$ the {\em proper Lorentzian group} 
consisting of all matrices of the form 
\[
A(\phi) = \left( \begin{array}{cc}
        \ch\phi & \sh\phi \\
	\sh\phi & \ch\phi
       \end{array}
     \right),
\]
where $\phi \in \R$,  see \cite{BN}, \cite{Y}.

In $L^2$ a Lorentzian vector $\underline{u}$ is called to 
{\em spacelike},
{\em lightlike} or {\em timelike} 
if $\langle\underline{u},\underline{u}\rangle_L >0$,
 $\langle\underline{u},\underline{u}\rangle_L =0$ or 
$\langle\underline{u},\underline{u}\rangle_L <0$, respectively.

We say that a timelike vector $\underline{u}=(u_1,u_2)$ is 
{\em future--pointing} or 
{\em past--pointing} if $u_2>0$ or $u_2<0$, respectively.
Similarly, a spacelike  vector $\underline{u}=(u_1,u_2)$ is 
{\em future--pointing} or 
{\em past--pointing} if $u_1>0$ or $u_1<0$, respectively.
The following Lemma was proven in \cite{BN}.

\begin{lemma} \label{revers_timelike} (Reversed triangle inequality) 
Let $\underline{x}$ and $\underline{y}$ be future--pointing timelike vectors in $L^2$.
Then $\underline{x}+\underline{y}$ is a future--pointing timelike vector and 
$$
\parallel \underline{x}+\underline{y} \parallel \geq  \parallel \underline{x} \parallel + \parallel \underline{y}\parallel .
$$
Here the equality holds iff $\underline{y}=c\underline{x}$ for some $c>0$.
\end{lemma}

\begin{lemma} \label{revers_spacelike}
 
Let $\underline{x}$ and $\underline{y}$ be future--pointing spacelike vectors in  $L^2$.
Then $\underline{x}+\underline{y}$ is a future--pointing spacelike vector and 
$$
\parallel \underline{x}+\underline{y} \parallel \geq  \parallel \underline{x} \parallel + \parallel \underline{y}\parallel .
$$
Here the equality holds iff $\underline{y}=c\underline{x}$ for some $c>0$.
\end{lemma}

\proof 
Suppose that $\underline{x}=\overrightarrow{AB}$, $\underline{y}=\overrightarrow{BC}$
 and $\underline{x}+\underline{y}=\overrightarrow{AC}$. Now after reflecting the points $A$, $B$ and $C$ 
to the line $\underline{y}=\underline{x}$ we get the points $A'$, $B'$ and $C'$. Let
 $\underline{d}:=\overrightarrow{A'B'}$ and  $\underline{e}:=\overrightarrow{B'C'}$. 
 Then $\underline{d}+\underline{e}=\overrightarrow{A'C'}$ and it comes from the definition of the norm that
$$
\parallel \underline{x} \parallel = \parallel \underline{d}\parallel,\ \parallel \underline{y} \parallel = \parallel \underline{e}\parallel\mbox{ and }\parallel \underline{d}+\underline{e} \parallel =\parallel \underline{x}+\underline{y} \parallel.
$$
Hence using Lemma \ref{revers_timelike} we get that
$$
\parallel \underline{d}+\underline{e} \parallel = \parallel \underline{x}+\underline{y} \parallel \geq \parallel \underline{x} \parallel + \parallel \underline{y}\parallel=\parallel \underline{d} \parallel + \parallel \underline{e}\parallel. 
$$
\qed

\begin{cor}
Suppose that $OABC$ is a quadrilaterial  where 
$a=\parallel \overrightarrow{OA}\parallel$, $b=\parallel \overrightarrow{BC}\parallel$, 
$g=\parallel \overrightarrow{OC}\parallel$, $h=\parallel \overrightarrow{AB}\parallel$. 
Suppose that $\overrightarrow{OA}$, $\overrightarrow{BC}$, $\overrightarrow{OC}$
and $\overrightarrow{AB}$ are future--pointing spacelike vectors 
(see Figure \ref{planar_4R_linkage}). Then $g\geq a+b+h$. 

Moreover, if $g=a+b+h$, then the vectors $\overrightarrow{OA}$, $\overrightarrow{BC}$,  
$\overrightarrow{OC}$ and $\overrightarrow{AB}$ are collinear.
\end{cor}
\proof 
The inequality $g\geq a+b+h$ is a direct consequence of Lemma \ref{revers_spacelike}. 

Now, assume that $g=a+b+h$. We prove first that the vectors $\overrightarrow{OC}$ and $\overrightarrow{OB}$
are collinear. Suppose, indirectly, that  $\overrightarrow{OC}$ and $\overrightarrow{OB}$
are not collinear. Since $\overrightarrow{OB}$ is a future--pointing spacelike
 vector, hence 
we can apply the reversed triangle inequality (Lemma \ref{revers_spacelike}) 
for the triangle $OBC$. Then it follows that
$$
g=\parallel \overrightarrow{OC} \parallel > \parallel \overrightarrow{OB} \parallel+ 
\parallel \overrightarrow{BC} \parallel = b+ \parallel \overrightarrow{OB} \parallel
$$
We can apply again Lemma \ref{revers_spacelike} for the triangle $OAB$. Hence
$$
\parallel \overrightarrow{OB} \parallel \geq a+h.
$$
Consequently 
$$
g=\parallel \overrightarrow{OC} \parallel > a+b+h,
$$
which is a contradiction. 

Similar argument shows that $\overrightarrow{OB}$ and $\overrightarrow{OA}$ are colllinear.
Hence $\overrightarrow{OC}$ and $\overrightarrow{OA}$ are collinear.
We can use very similar arguments to prove that $\overrightarrow{OA}$ and 
$\overrightarrow{AB}$ are collinear vectors. So the result follows. \qed

We define now the notion of angle on the Minkowskian plane. 

Let $\underline{x}$ and $\underline{y}$ be two future--pointing timelike unit vectors.
We say that $\phi\in \R$ is the {\em (oriented) angle}
 from  $\underline{x}$ to $\underline{y}$ if 
$A(\phi)\underline{x}=\underline{y}$. The {(unoriented) angle between
$\underline{x}$ and $\underline{y}$} is defined to be $|\phi|$.
Then it comes from the definition that
$$
\ch\phi=- \langle\underline{x}, \underline{y}\rangle_L,
$$
where the right-hand side is greater than $1$.

When $\underline{x}$ and $\underline{y}$ are future--pointing 
timelike vectors, than the angle $\phi$ for 
$\underline{x}$ and $\underline{y}$
is the same as for $\underline{x}/\parallel \underline{x}\parallel$ and
 $\underline{y}/\parallel \underline{y}\parallel$. 

We have 
$$
\ch\phi=- \frac{\langle\underline{x}, \underline{y}\rangle_L}{\parallel \underline{x}\parallel \parallel \underline{y}\parallel}
$$

We can give the same definition for the angle between future--pointing spacelike vectors. 
From this definition, we obtain similar formula for $\ch\phi$: 
if $\underline{x}$ and $\underline{y}$ are two future--pointing spacelike unit vectors,
then   
$$
\ch\phi= \langle\underline{x}, \underline{y}\rangle_L,
$$
where $\phi$ is the (oriented) angle from  $\underline{x}$ to $\underline{y}$.  
%The angle between two timelike lines (though the zero vector) is the angle between 
%the two future--pointing timelike unit vectors which lie on the respective lines.

In the Euclidean plane, a motion can be represented by a combination of 
a rotation and translation. It is well--known that any motion can be expressed
using the matrix operation
$$
M(t,a,b)\left( \begin{array}{c}
                x_1 \\ 
                x_2 \\
                1 \\
               \end{array}\right) =
\left( \begin{array}{ccc}
        \mbox{cos}(t) & -\mbox{sin}(t) & a \\
	\mbox{sin}(t) & \mbox{cos}(t) & b \\
                 0    &        0      & 1 \\ 
    \end{array}
     \right)
\left( \begin{array}{c} x_1 \\ 
                x_2 \\
                1 \\
               \end{array}\right)
$$

Similarly, for the Minkowskian plane the group of motions is 
the following:
$$
\overline{G}:=\left\{ \left( \begin{array}{ccc}
        \ch t & \sh t & a \\
	\sh t & \ch t & b \\
         0    &  0    & 1 \\ 
    \end{array}
     \right):~ t,a,b\in\R \right\}
$$

We now introduce an important class of triangles.
By a {\em pure triangle} we mean a triangle 
with vertices $A, B$ and $C$ such that
$\overrightarrow{AB}$ and $\overrightarrow{BC}$ 
are future--pointing timelike vectors. In the following we assume 
that we named the vertices of a pure triangle $ABC$ in this 
manner. The angle 
%$\widehat{A}$ is the angle between the lines 
%$AB$ and $AC$ and the angle
$\widehat{C}$ is the angle between
the lines $BC$ and $AC$. 

Finally we recall here for the Minkowskian cosine rule:
\begin{thm}(see \cite[Theorem 7]{BN}) \label{Mink_cos_law} 
 Let $\triangle ABC$ be a pure triangle. Then 
\begin{equation}
c^2=a^2+b^2-2ab {\rm ch}(\widehat{C})
\end{equation}
where $a=\parallel \overrightarrow{BC}\parallel$, 
$b=\parallel \overrightarrow{AC}\parallel$ and 
$c=\parallel \overrightarrow{AB}\parallel$.
\end{thm}
\qed

\section{The main results}
\subsection{Position Analysis of the Minkowskian 4R linkage}

Recall that the four-bar linkage is a mechanism that lies in the plane 
and consists of four bars 
connected by joints that allow rotation only in the plane of the mechanism, 
see Figure \ref{planar_4R_linkage}.

Throughout this Chapter we suppose that $\overrightarrow{OA}$, 
 $\overrightarrow{OC}$ and $\overrightarrow{CB}$
are future--pointing spacelike vectors. Suppose that $\overrightarrow{AB}$ 
is a spacelike vector.

Let the fixed and the moving pivots of 
the input crank be $O$ and $A$, respectively.

Let the fixed and the moving pivots of the output crank be $C$ and $B$, 
respectively. We define 
the distances between these points as follows:
$$
a:=\parallel \overrightarrow{OA}\parallel,\ 
b:=\parallel \overrightarrow{BC}\parallel ,\ 
g:=\parallel \overrightarrow{OC}\parallel,\
h:=\parallel \overrightarrow{AB}\parallel.
$$

To analyse the linkage, we locate the origin in the fixed Minkowskian 
frame $F$ at $O$ and 
orient it so that the $x$-axis passes through the other fixed pivot $C$.

\begin{thm} \label{psi}
Let $\theta$ be the input angle measured around $O$ from the $x$-axis of $F$ to
$OA$. Let $\psi$ be the angular position of the output crank $CB$ (see Figure \ref{planar_4R_linkage}). Then
\begin{equation} \label{psi2}
\psi=2\mbox{artanh}\frac{-B(\theta)\pm \sqrt{B(\theta)^2+C(\theta)^2-A(\theta)^2}}{A(\theta)+C(\theta)}, 
\end{equation}
where
$$
A(\theta)=2gb-2ab\ch\theta,
$$
$$
B(\theta)=2ab\sh\theta,
$$
and
$$
C(\theta)=h^2-g^2-b^2-a^2+2ag\ch\theta.
$$
\end{thm}

\proof
Since $h=\parallel \overrightarrow{AB}{\parallel}$ is constant, 
we get the constraint equation as
\begin{equation} \label{main}
\parallel \overrightarrow{AB}{\parallel}^2=h^2. 
\end{equation}
 
It is easy to verify that the coordinates of
 $A$ and $B$ is given by $A=(a \ch\theta, a\sh\theta)$ and
\begin{equation} \label{B_eq}
B=(g+b\ch\psi, b\sh\psi).
\end{equation}

If we substitute these coordinates into (\ref{main}), then we obtain
$$
b^2+g^2+a^2+2gb\ch\psi -2ag\ch\theta -2ab\ch\psi \ch\theta
+2ab\sh\psi \sh\theta=h^2.
$$

If we gather the coefficients of $\ch\psi$ and  $\sh\psi$, we obtain the constraint 
equation for the 4R chain as
\begin{equation} \label{main2}
A(\theta)\ch\psi +B(\theta)\sh\psi =C(\theta),
\end{equation}
where
$$
A(\theta)=2gb-2ab\ch\theta,
$$
$$
B(\theta)=2ab\sh\theta,
$$
and
$$
C(\theta)=h^2-g^2-b^2-a^2+2ag\ch\theta.
$$

We solve this equation using the tan-half-technique.
This technique uses a transformation of variables 
to convert $\ch\psi$ and $\sh\psi$ into algebraic functions
 of $\mbox{th}\left(\psi/2 \right)$.

Introduce the parameter $y=\mbox{th}\left(\psi/2 \right)$, 
which allows us to define
$$
\ch\psi=\frac{1+y^2}{1-y^2}\mbox{ and }\sh\psi=\frac{2y}{1-y^2}
$$ 
Substitute into (\ref{main2}) to obtain
$$
(A(\theta)+C(\theta))y^2+2B(\theta)y+A(\theta)-C(\theta)=0.
$$

This equation is solved using the quadratic formula to obtain
$$
\mbox{th}\left(\psi/2 \right)=\frac{-B(\theta)\pm \sqrt{B(\theta)^2+C(\theta)^2-A(\theta)^2}}{A(\theta)+C(\theta)}.
$$
\qed
\begin{remark}
It can be derived the following alternative equation for $\psi$:
$$
\psi=-\mbox{artanh}\frac{B(\theta)}{A(\theta)}\pm \mbox{arch}\frac{C(\theta)}{\sqrt{A(\theta)^2-B(\theta)^2}}.
$$

Namely we infer from (\ref{main2}) that 
\begin{equation} \label{main25}
\frac{A(\theta)}{\sqrt{A(\theta)^2-B(\theta)^2}}\ch\psi +\frac{B(\theta)}{\sqrt{A(\theta)^2-B(\theta)^2}}\sh\psi =\frac{C(\theta)}{\sqrt{A(\theta)^2-B(\theta)^2}}.
\end{equation} 
Since 
$$
\frac{A(\theta)}{\sqrt{A(\theta)^2-B(\theta)^2}} \geq 1, 
$$
hence there exists $\delta$ such that
$$
\ch\delta =\frac{A(\theta)}{\sqrt{A(\theta)^2-B(\theta)^2}}
$$
and 
$$
\sh\delta = \frac{B(\theta)}{\sqrt{A(\theta)^2-B(\theta)^2}}.
$$
Then clearly
$$
\tah\delta = \frac{B(\theta)}{A(\theta)}.
$$

We infer from (\ref{main25}) that
$$
\mbox{ch}(\delta + \psi)=\ch\delta \ch\psi + \sh\delta \sh\psi=\frac{C(\theta)}{\sqrt{A(\theta)^2-B(\theta)^2}}.
$$

The result follows. 
\end{remark}
\begin{remark}
A graph of $\psi$ as a function of $\theta$ is called 
the (kinematic) {\em transmission function}. The transmission function has a very 
characteristic shape, hence the type of the Minkowskian four-bar linkage 
can immediately be read from it.
\end{remark}

\begin{remark}
Notice that there are two angles $\psi$ for each angle $\theta$.
This arises because the moving pivot $B$ of the output crank can be assembled above or
below the diagonal joining the moving pivot $A$ of the input crank to the
fixed pivot of $C$.
\end{remark}

\subsection{Branching points}

The transmission function gives for any $\theta$ two, one, no or an infinite number 
of values of $\psi$. No values for $\psi$ exist iff
 the position of the linkage is impossible.
In an 'undetermined position' $\theta$ and $\psi$ are independent of each other: at a 
single $\theta$, an infinite amount of values of $\psi$ is possible or vice versa.
We call these $\theta$ angles to the {\em branching points} of the mechanism.
 
\begin{cor} (Branching points) \label{branching}
Suppose that $g\neq b$. Then  
the branching points of 
the Minkowskian planar 4R occours at 
\begin{equation} \label{branch_theta}
\theta = \mbox{arch}\left( \frac{a^2-h^2}{2a(g-b)}+\frac{g-b}{2a} \right).
\end{equation}
\end{cor}
\proof
It is clear from Theorem \ref{psi} that $\theta$ is a branching point of the mechanism iff
$$
A(\theta)+C(\theta)=0.
$$
This means that
$$
2gb-2ab\ch\theta +h^2-g^2-b^2-a^2+2ag\ch\theta=0,
$$
that is
$$
\ch\theta = \frac{a^2-h^2}{2a(g-b)}+\frac{g-b}{2a}.
$$
\qed
\begin{cor}  \label{branching2}
Suppose that $g=b$ and $h\neq a$. Then 
there is no branching points. 
\end{cor}

\proof
If $g=b$ and $h\neq a$, then 
$$
A(\theta)+C(\theta)=h^2-a^2\neq 0.
$$
for each $\theta$.
\qed
\begin{cor}  \label{branching3}
Suppose that $g=b$ and $h=a$. Then 
all points are branching points. 
\end{cor}
\proof 

This is clear, since then 
$$
A(\theta)+C(\theta)=0
$$
for each $\theta$.
\qed
\begin{remark}
The deeper investigation of branching points will be one of the topics
of our next article.   
\end{remark}

\subsection{Coupler Angles}

\begin{figure}[htp]
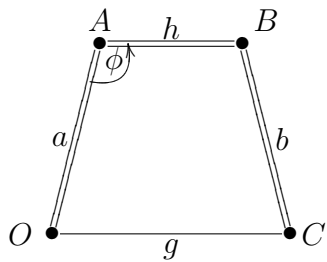

    \centering 
% coupler angle
\[
\xy
0; <3pt,0pt>:
(0,0)*{\bullet}="O"; 
(30,0)*{\bullet}="C";
(6,24)*{\bullet}="A"; 
(24,24)*{\bullet}="B";
(40,0)*{}="Y";
(35,0)*{}="K";
(21,36)*{}="T";
(22.5,30)*{}="S";
"O"+(-4,0)*{O}; 
"C"+(3,0)*{C};
"A"+(0,3)*{A};
"B"+(3,3)*{B};
"O"; "A"; **\dir{=}; ?(.8)*{}="X2"; ?(.2)*{}="X3";
"O"; "C"; **\dir{-}; ?(.8)*{}="X1";  ?(.2)*{}="X8";
"A"; "B"; **\dir{=}; ?(.8)*{}="X4"; ?(.3)*{}="X5"; 
"B"; "C"; **\dir{=};  ?(.8)*{}="X6"; ?(.2)*{}="X7";
%"C"; "Y"; **\dir{-};
%"B"; "T"; **\dir{-};
%"X1"; "X2"; **\crv{(5,5)};
%"X3"; "X4"; **\crv{(11,19)};
{\ar@/_0.6pc/ "X3";"X4"};
%"X5"; "X6"; **\crv{(19,19)};
%"X7"; "X8"; **\crv{(25,5)};
%"X5"; "S"; **\crv{(17,28)};
%"X7"; "K"; **\crv{(34,5)};
%"O"+(3,2)*{\theta};
%"C"+(2,2)*{\psi};
%"B"+(-3,3)*{\zeta};
"A"+(1.8,-2.2)*{\phi};
(15,-2)*{g};
(1,12)*{a};
(15,26)*{h};
(29,12)*{b};
\endxy
\]
\caption[The coupler angle]{The coupler angle}
\label{coupler_angle_phi}
\end{figure}

Let $\phi$ denote the angle between the vectors $ \overrightarrow{AO}$ and $ \overrightarrow{AB}$
 (see Figure \ref{coupler_angle_phi}).

\begin{thm} \label{coupler}
If $\phi$ denotes the coupler angle, then
$$
\phi=\mbox{artanh}\left( \frac{b\sh\psi- a\sh\theta}{g+b\ch\psi-a\ch\theta} \right)-\theta + \pi
$$
\end{thm}

\proof

Since $\phi$ was the coupler angle, hence $\theta +\phi -\pi$ is the angle to $AB$ 
from the $x$-axis of $F$. Consequently we can write the coordinates of $B$ in terms of
$\phi$ as
\begin{equation} \label{B_eq2}
B=(a\ch\theta+ h\mbox{ch}\left( \theta+\phi-\pi \right), a\sh\theta+ h\mbox{sh}\left( \theta+\phi-\pi \right) ).
\end{equation} 

If we equate the two forms (\ref{B_eq}) and (\ref{B_eq2}) for $B$, we obtain the following 
loop equations of the four-bar linkage:
$$
a\ch\theta+ h\mbox{ch}\left( \theta+\phi-\pi \right)=g+b\ch\psi
$$
$$
a\sh\theta+ h\mbox{sh}\left( \theta+\phi-\pi \right)=b\sh\psi
$$

Therefore, for a given value of the driving crank $\theta$, $\mbox{ch}\left( \theta+\phi-\pi \right)$
and $\mbox{sh}\left( \theta+\phi-\pi \right)$ are given by
$$
\mbox{ch}\left( \theta+\phi-\pi \right)=\frac{g+b\ch\psi-a\ch\theta}{h}
$$
and
$$
\mbox{sh}\left( \theta+\phi-\pi \right)=\frac{b\sh\psi- a\sh\theta}{h}.
$$
Hence 
$$
\mbox{th}\left( \theta+\phi-\pi \right)=\frac{b\sh\psi- a\sh\theta}{g+b\ch\psi-a\ch\theta}
$$
and Theorem \ref{coupler} follows.
\qed
\begin{remark}
Notice that we obtain a unique value for $\phi$ associated 
with each of the two 
solutions for the output angle $\psi$.
\end{remark}

\subsection{The transmission angle}

\begin{figure}[htp]
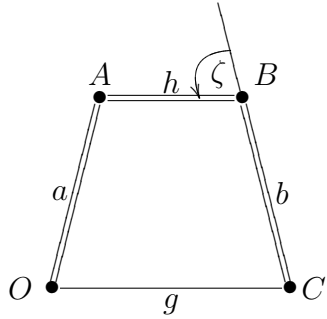

    \centering 
\[
\xy
0; <3pt,0pt>:
(0,0)*{\bullet}="O"; 
(30,0)*{\bullet}="C";
(6,24)*{\bullet}="A"; 
(24,24)*{\bullet}="B";
(40,0)*{}="Y";
(35,0)*{}="K";
(21,36)*{}="T";
(22.5,30)*{}="S";
"O"+(-4,0)*{O}; 
"C"+(3,0)*{C};
"A"+(0,3)*{A};
"B"+(3,3)*{B};
"O"; "A"; **\dir{=}; ?(.8)*{}="X2"; ?(.2)*{}="X3";
"O"; "C"; **\dir{-}; ?(.8)*{}="X1";  ?(.2)*{}="X8";
"A"; "B"; **\dir{=}; ?(.8)*{}="X4"; ?(.3)*{}="X5"; 
"B"; "C"; **\dir{=};  ?(.8)*{}="X6"; ?(.2)*{}="X7";
%"C"; "Y"; **\dir{-};
"B"; "T"; **\dir{-};
%"X1"; "X2"; **\crv{(5,5)};
%"X3"; "X4"; **\crv{(11,19)};
%"X5"; "X6"; **\crv{(19,19)};
%"X7"; "X8"; **\crv{(25,5)};
%"X5"; "S"; **\crv{(17,28)};
{\ar@/_0.6pc/ "S";"X5"};
%"X7"; "K"; **\crv{(34,5)};
%"O"+(3,2)*{\theta};
%"C"+(2,2)*{\psi};
"B"+(-3,3)*{\zeta};
%"A"+(2,-2)*{\phi};
(15,-2)*{g};
(1,12)*{a};
(15,26)*{h};
(29,12)*{b};
\endxy
\]
\caption[The transmission angle]{The transmission angle}
\label{transmission_angle_phi}
\end{figure}

Let $\zeta$ denote the angle between the coupler and the driven crank 
at $B$, the {\em transmission angle} of the linkage (see Figure \ref{transmission_angle_phi}).

\begin{thm} \label{transmission}
If $\zeta$ is the transmission angle of the linkage, then
$$
\zeta =\mbox{arch}\left( \frac{-g^2-a^2+h^2+b^2+2ag\ch\theta}{2bh}\right).
$$
\end{thm}

\proof 

Let $d:=\parallel \overrightarrow{AC} \parallel$.

To determine $\zeta$ in terms of $\theta$, equate the Minkowskian cosine laws
(Theorem \ref{Mink_cos_law}) for the diagonal $AC$ for the triangles $\bigtriangleup COA$ and
$\bigtriangleup ABC$. Since $\zeta$ is the interior angle at $B$, we have
$$
d^2=g^2+a^2-2ag\ch\theta =h^2+b^2-2bh\ch\zeta.
$$
The result follows
\qed

\subsection{The coupler curve}

\begin{figure}[htp]
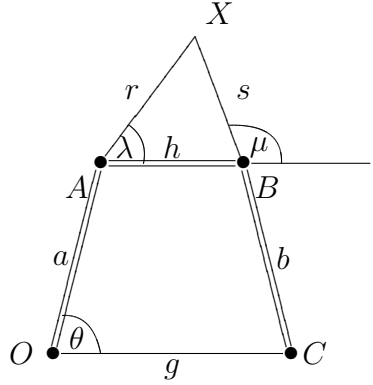

    \centering 
% coupler curve
\[
\xy
0; <3pt,0pt>:
(0,0)*{\bullet}="O"; 
(30,0)*{\bullet}="C";
(6,24)*{\bullet}="A"; 
(24,24)*{\bullet}="B";
(40,0)*{}="Y";
(35,0)*{}="K";
(21,36)*{}="T";
(22.5,30)*{}="S";
(18,40)*{}="X";
(40,24)*{}="W";
"O"+(-4,0)*{O}; 
"C"+(3,0)*{C};
"A"+(-3,-3)*{A};
"B"+(3,-3)*{B};
"X"+(3,3)*{X};
"O"; "A"; **\dir{=}; ?(.8)*{}="X2"; ?(.2)*{}="X3";
"O"; "C"; **\dir{-}; ?(.8)*{}="X1";  ?(.2)*{}="X8";
"A"; "B"; **\dir{=}; ?(.7)*{}="X4"; ?(.3)*{}="X5"; 
"B"; "C"; **\dir{=};  ?(.8)*{}="X6"; ?(.2)*{}="X7";
%"C"; "Y"; **\dir{-};
"A"; "X"; **\dir{-}; ?(.7)*{}="X9";
"B"; "X"; **\dir{-}; ?(.7)*{}="G";
"B"; "W"; **\dir{-};  ?(.7)*{}="F";
"X1"; "X2"; **\crv{(5,5)};
"F"; "G"; **\crv{"B"+(5,5)};
%"X5"; "X6"; **\crv{(19,19)};
%"X7"; "X8"; **\crv{(25,5)};
%"X5"; "S"; **\crv{(17,28)};
%"X7"; "K"; **\crv{(34,5)};
"X9"; "X4"; **\crv{(12,27)};
"O"+(3,2)*{\theta};
%"C"+(2,2)*{\psi};
%"B"+(-3,3)*{\zeta};
"A"+(3,2)*{\lambda};
"B"+(2,2)*{\mu};
"A"+(4,9)*{r};
"A"+(18,9)*{s};
(15,-2)*{g};
(1,12)*{a};
(15,26)*{h};
(29,12)*{b};
\endxy
\]
\caption[The coupler curve]{The coupler curve}
\label{coupler_curve2}
\end{figure}

In this subsection, we study the motion of the coupler by analyzing 
the curve traced by a point on the coupler link. We get the parametrized equation 
of this curve, the {\em coupler curve} from the kinematics equations
of the driving RR chain (see Figure \ref{coupler_curve2}).

\begin{thm} \label{coupler_curve} 
The curve traced by any point of the coupler link 
of a Minkowskian planar four-bar
linkage is algebraic, of sixth degree.
\end{thm}

\proof

Let $\underline{x}=(x,y)^T$ be the coordinates of a coupler point 
in the frame $M$ located at $A$ with its $x$-axis along $AB$.

Let $\underline{X}=(X,Y)^T$ be the coordinates of a coupler point 
in the frame $F$. 
%located at $O$ with its $x$-axis along $OC$.

We obtain the algebraic equation of the coupler curve 
by defining the coordinates of $\underline{X}=(X,Y)^T$ from two points of
view. Let the coupler triangle $\bigtriangleup XAB$ have length $r$ and $s$
given by 
$$
r=\parallel \overrightarrow{AX} \parallel=\sqrt{x^2-y^2}
$$ 
and 
$$
s=\parallel \overrightarrow{BX} \parallel=\sqrt{(x-h)^2-y^2}.
$$
If $\lambda$ is the angle to $AX$ in $F$ and $\mu$ is the angle to
$BX$ (see Figure \ref{coupler_curve2}), then we have by definition that
$$
\overrightarrow{AX}=(r\ch \lambda,r\sh\lambda)
$$
and
$$
\overrightarrow{BX}=(s\ch \mu,s \sh\mu).
$$
If we substitute into the identities 
$\parallel \overrightarrow{OA} \parallel^2=a^2$ and 
$\parallel \overrightarrow{CB} \parallel^2=b^2$ and rearrange these equations, we obtain
\begin{equation} \label{coupler_equation1}
(X-s\ch\mu -g)^2-(Y-s\sh\mu)^2=h^2
\end{equation}
and
\begin{equation} \label{coupler_equation2}
(X-r\ch\lambda)^2-(Y-r\sh\lambda)^2=a^2.
\end{equation}

If we expand equations (\ref{coupler_equation1}) and (\ref{coupler_equation2}), we get
\begin{equation} \label{coupler_equation3}
2Xs\ch\mu -2Ys\sh\mu -2gs\ch\mu -X^2+2Xg+Y^2-g^2+h^2-s^2=0 
\end{equation}
and 
\begin{equation} \label{coupler_equation4}
2Xr\ch\lambda-2Yr\sh\lambda-X^2+Y^2+a^2-r^2=0.
\end{equation}
Let  $\gamma:=\mu-\lambda$ and substitute $\mu=\lambda+\gamma$ 
into the equation (\ref{coupler_equation3}).

If we rearrange these equations we get
$$
A_1(X,Y)\ch\lambda +B_1(X,Y)\sh\lambda=C_1(X,Y),
$$
\begin{equation} \label{Cupler_curve_eq_sys}
A_2(X,Y)\ch\lambda +B_2(X,Y)\sh\lambda=C_2(X,Y),
\end{equation}
where
$$
A_1(X,Y):=2Xs\sh\gamma-2Ys\ch\gamma-2gs\sh\gamma,
$$
$$
B_1(X,Y):=2Xs\ch\gamma-2Ys\sh\gamma-2gs\ch\gamma,
$$
$$
C_1(X,Y):=X^2-2Xg-Y^2+g^2-h^2+s^2,
$$
$$
A_2(X,Y):=2rX,
$$
$$
B_2(X,Y):=-2rY,
$$
and
$$
C_2(X,Y):=X^2-Y^2-a^2+r^2.
$$

We can eliminate $\lambda$ from the equation system (\ref{Cupler_curve_eq_sys}) 
by solving linearly for $u=\ch\lambda$ 
and $v=\sh\lambda$. Then we can impose the condition
$u^2-v^2=1$. This yields to
$$
(C_1(X,Y)B_2(X,Y)-C_2(X,Y)B_1(X,Y))^2-(A_2(X,Y)C_1(X,Y)-A_1(X,Y)C_2(X,Y))^2-
$$
$$
-(A_1(X,Y)B_2(X,Y)-A_2(X,Y)B_1(X,Y))^2=0.
$$
Notice that $A_i(X,Y)$ and $B_i(X,Y)$ are linear in the coordinates $X$ and $Y$,
 and $C_i(X,Y)$ are quadratic. Therefore the equation defines a curve of degree
6.
\qed
\subsection{The input and the output Crank Angles and the classification of the 4R
linkage}

The following Theorem describes the upper and the lower 
limiting angles.
\begin{thm} \label{input_Crank_Angles}
The upper and the lower 
limiting angles $\theta_{max}$ and $\theta_{min}$ that
define the range of movement of the input crank,
$$
\mbox{ch}(\theta_{min})=\frac{a^2+g^2-(b+h)^2}{2ag}
$$
and 
$$
\mbox{ch}(\theta_{max})=\frac{a^2+g^2-(b-h)^2}{2ag}.
$$
\end{thm}

\proof
The formula (see 
equation (\ref{psi2}) in Theorem \ref{psi}) 
that defines the output angle $\psi$ for a given input angle $\theta$ 
has a real solution if and only if $B(\theta)^2+C(\theta)^2-A(\theta)^2\geq 0$.
We obtain the maximum and the minimum
values for $\theta$ if we
set this condition to zero, which yields 
the following quadratic equation in $\ch\theta$
$$
4a^2g^2\mbox{ch}^2 (\theta)-4ag(g^2+a^2-h^2-b^2)\ch\theta +
$$ 
$$
+((g^2+a^2)-(h+b)^2)((g^2+a^2)-(h-b)^2)=0.
$$
The roots of this equation are the given equations for $\mbox{ch}(\theta_{min})$ 
and $\mbox{ch}(\theta_{max})$.
\qed
\begin{remark}
These equations are the cosine laws for the two ways that the triangle 
$AOC\bigtriangleup$ can be formed with the coupler $AB$ aligned with the output
crank $CB$, see Figure \ref{fig_coslaw}, Figure \ref{fig_coslaw2} and Figure \ref{fig_coslaw3}.
\end{remark}

\begin{figure}[htp] 
    \centering 
% theta_max
\[
\xy
0; <3pt,0pt>:
(0,0)*{\bullet}="O"; 
(30,24)*{\bullet}="A"; 
(40,0)*{\bullet}="C"; 
"O"+(-4,0)*{O}; 
"C"+(3,0)*{C};
"A"+(0,3)*{A};
"O"; "C"; **\dir{-}; ?(.65)*{}="F";
"O"; "A"; **\dir{=}; ?(.65)*{}="G";
"A"; "C"; **\dir{=}; ?(.4)*{\bullet}="B";
"F"; "G"; **\crv{(14,6)};
"B"+(2,3)*{B};
"O"+(9,2)*{\theta_{max}};
(20,-2)*{g};
(11,12)*{a};
(36,18)*{h};
(41,7)*{b};
\endxy
\]
\caption[The_planar_4R_linkage2]{The angles $\theta_{min}$ and $\theta_{max}$
are the limits to the range of movement of the input link}
\label{fig_coslaw}
\end{figure} 

\begin{figure}[htp]
    \centering 
\[
\xy
0; <3pt,0pt>:
(0,0)*{\bullet}="O"; 
(30,24)*{\bullet}="A"; 
(25,36)*{\bullet}="B";
(40,0)*{\bullet}="C"; 
"O"+(-4,0)*{O}; 
"C"+(3,0)*{C};
"A"+(2,3)*{A};
"B"+(0,3)*{B};
"O"; "C"; **\dir{-}; ?(.65)*{}="F";
"O"; "A"; **\dir{=}; ?(.65)*{}="G";
"B"; "C"; **\dir{=};
"F"; "G"; **\crv{(15,6)};
"O"+(9,2)*{\theta_{min}};
(20,-2)*{g};
(11,12)*{a};
(31,32)*{h};
(42,13)*{b-h};
\endxy
\]
\caption[The planar_4R_linkage3]{The angles $\theta_{min}$ and $\theta_{max}$
are the limits to the range of movement of the input link}
\label{fig_coslaw2}
\end{figure} 

\begin{figure}[htp] 
    \centering 
\[
\xy
0; <3pt,0pt>:
(0,0)*{\bullet}="O"; 
(30,24)*{\bullet}="A"; 
(45,-12)*{\bullet}="B";
(40,0)*{\bullet}="C"; 
"O"+(-4,0)*{O}; 
"C"+(4,0)*{C};
"A"+(2,3)*{A};
"B"+(3,0)*{B};
"O"; "C"; **\dir{-}; ?(.65)*{}="F";
"O"; "A"; **\dir{=}; ?(.65)*{}="G";
"A"; "B"; **\dir{=};
"F"; "G"; **\crv{(15,6)};
"O"+(9,2)*{\theta_{min}};
(20,-2)*{g};
(11,12)*{a};
(47,-6)*{b};
(42,13)*{h-b};
\endxy
\]
\caption[The planar_4R_linkage3]{The angles $\theta_{min}$ and $\theta_{max}$
are the limits to the range of movement of the input link}
\label{fig_coslaw3}
\end{figure}

The hyperbolic cosine function does not distinguish between 
$\pm \theta$, so there are actually two limits for 
%$\ch\theta_{min}$ and 
$\mbox{ch}(\theta_{max})$ above and below $OC$. 

If $\theta_{max}$ does not exist, then the crank has no lower limit 
to its movement and it rotated though $\theta=0$ to reach negative 
values below the segment $OC$.  

Thus $\mbox{ch}(\theta_{max})<1$ is the condition that there is no lower limit 
to the input crank rotation, that is, 
$$
\frac{a^2+g^2-(b-h)^2}{2ag}<1.
$$
We can simplify this to yield
$$
(b-h)^2-(a-g)^2>0,
$$
If we factor the difference to two squares, we obtain
$$
(-a+g+b-h)(a-g+b-h)>0,
$$
that is
$$
T_1T_2>0,
$$
where $T_1:=-a+g+b-h$ and $T_2:=a-g+b-h$. 

If $\mbox{ch}(\theta_{min})>1$, then both upper and lower 
limiting angles exist. The input crank rocks in one of three separate ranges:
(i) $\theta_{max} \leq\ \theta\leq \infty$, 
(ii) $-\theta_{min} \leq\theta \leq \theta_{min}$
or (iii) $-\infty \leq\theta\leq -\theta_{max}$.

It is easy to see that the $\mbox{ch}(\theta_{min})>1$ condition is equivalent with 
$$
(g-a-b-h)(g-a+b+h)>0,
$$ 
that is,
$T_3T_4>0$, where $T_3:=g-a-b-h$ and $T_4:=g-a+b+h$. 

Finally we get the following result:

\begin{thm} \label{input_classification}
We can identify three types of movement available to the input crank of a 4R 
linkage:

1. A {\bf crank}: $T_1T_2\geq 0$ and $T_3T_4\leq 0$, in which case no $\theta_{min}$, $\theta_{max}$ exists, and the input 
crank rocks though $\theta=0$.
%between the values $\pm \theta_{max}$. 

2. A {\bf rocker}: $T_1T_2<0$ and $T_3T_4\leq 0$, which means that both upper and lower 
limiting angles exist, and the crank cannot pass though $0$. Instead, it rocks 
in one of the two separate ranges: (i) $\theta_{max} \leq\theta \leq \infty$
or (ii) $-\infty \leq\theta \leq -\theta_{max}$.

3. A {\bf superrocker}: $T_1T_2<0$ and $T_3T_4>0$, which means that both upper and lower 
limiting angles exist. It rocks in one of the three separate ranges: 
(i) $\theta_{max} \leq\ theta\leq \infty$, 
(ii) $-\theta_{min} \leq\theta \leq \theta_{min}$
or (iii) $-\infty \leq\theta\leq -\theta_{max}$.  
\end{thm}
\qed

The range of movement of the output crank can be analyzed in the same way.
The limiting positions occur when the input crank $OA$ and the coupler $AB$ 
become aligned.

\begin{thm} \label{output_Crank_Angles}
The limits  $\psi_{max}$ and $\psi_{min}$ are defined 
by the equations
$$
\mbox{ch}(\psi_{max})=\frac{(a+h)^2-g^2-b^2}{2bg}
$$
and
$$
\mbox{ch}(\psi_{min})=\frac{(a-h)^2-g^2-b^2}{2bg}.
$$
\end{thm}
\qed

\begin{remark}
In this case $\psi$ is the exterior angle and changes the 
sign of the hyperbolic cosine term in the cosine law formula.
\end{remark}

We find the condition for no lower limit $\psi_{max}$ is
$$
\mbox{ch}(\psi_{max})=\frac{(a+h)^2-g^2-b^2}{2bg}<1.
$$
Hence 
$$
(a+h)^2<(b+g)^2
$$
that is
$$
(g+b+a+h)(g+b-h-a)>0.
$$
Using the notations  
$$
T_1:=g+b-h-a\mbox{ and }T_7:=g+b+a+h,
$$
we get
$$
T_1T_7>0.
$$

But clearly $T_7>0$. 

It is easy to see that the condition 
$\mbox{ch}(\psi_{min})>1$
is equivalent with $(a-h-g-b)(a-h+g+b)>0$, that is $T_4T_5<0$, where 
$T_4:=-a+h+g+b$ and $T_5:=a-h+g+b$. 

Finally we get the following result.

\begin{thm} \label{output_classification}
We can identify three types of movement available to the output crank of a 4R 
Minkowskian linkage:

1. A {\bf crank}: $T_1\geq 0$ and $T_4T_5\geq 0$, in which case no $\psi_{min}$, $\psi_{max}$ exists, and the output 
crank rocks though $\psi=0$. 
%between the values $\pm \psi_{max}$. 

2. A {\bf rocker}: $T_1<0$ and $T_4T_5\geq 0$, which means that both upper and lower 
limiting angles exist, and the crank cannot pass though $0$. Instead, it rocks 
in one of the two separate ranges: (i) $\psi_{max} \leq\psi \leq \infty$
or (ii) $-\infty \leq\psi \leq -\psi_{max}$.

3. A {\bf superrocker}: $T_1<0$ and $T_4T_5<0$, which means that both upper and lower 
limiting angles exist. It rocks in one of the three separate ranges: 
(i) $\psi_{max} \leq\ \psi\leq \infty$, 
(ii) $-\psi_{min} \leq\psi \leq \psi_{min}$
or (iii) $-\infty \leq\psi\leq -\psi_{max}$.  

\end{thm}
\qed

Hence we can classify a Minkowskian planar 4R linkage by the movement of 
input and output cranks.
For example, a crank-rocker has a rotatable input link though $O$ and an output link 
that rocks.

Now we define special subclasses of the Minkowskian planar 4R mechanisms.

We say that a planar 4R is {\em normal}, if all the conditions $a+b+h\leq g$, $a+b+g\leq h$, $a+g+h\leq b$
and $b+g+h\leq a$ are satisfied. We say that a planar 4R is {\em strange} if it is 
not normal. This means that there exists one link, 
which is longer than the sum of other links' lengths. 

We say that a normal planar 4R is {\em rigid}, if one of the following 
conditions is satisfied: $a+b+h=g$, $a+b+g= h$, $a+g+h= b$
and $b+g+h= a$. 

We say that a nonrigid normal planar 4R is {\em irreducible}, if 
$T_1\neq 0$ and $T_2\neq 0$, i.e., $a+h\neq b+g$ and $a+b\neq g+h$.

Finally we say that a nonrigid normal planar 4R, which is not irreducible, is a 
{\em reducible} planar 4R.

First we investigate the strange mechanisms.
\begin{thm} Suppose that $g>a+b+h$. Then 
the Minkowskian planar 4R has a superrocker--crank type. 
\end{thm}
\proof

Since $T_1>0$, $T_2<0$, $T_3>0$, $T_4>0$ and $T_5>0$.

\qed
\begin{thm} Suppose that $a>g+b+h$. Then 
the Minkowskian planar 4R has a superrocker--superrocker type. 
\end{thm}
\proof

Since $T_1<0$, $T_2>0$, $T_3<0$, $T_4<0$ and $T_5>0$.

\qed
\begin{thm} Suppose that $h>g+b+a$. Then 
the Minkowskian planar 4R has a crank--superrocker type. 
\end{thm}
\proof

Since $T_1<0$, $T_2<0$, $T_3<0$, $T_4>0$ and $T_5<0$.

\qed
\begin{thm} Suppose that $b>g+h+a$. Then 
the Minkowskian planar 4R has a crank--crank type. 
\end{thm}
\proof

Since $T_1>0$, $T_2>0$, $T_3<0$, $T_4>0$ and $T_5>0$.
\qed

We summarized the classification of strange mechanisms in Table \ref{Strange}.

\begin{table}[htp] 
\begin{center}
\begin{tabular}{||c|l|l||}
\hline
 & Linkage Type & Condition   \\
\hline
1 & Crank--Crank & $b>g+h+a$  \\
2 & Crank--Superrocker & $h>g+b+a$   \\
3 & Superrocker--crank & $g>a+b+h$ \\
4 & Superrocker--superrocker & $a>g+b+h$ \\

\hline
\end{tabular}
\vspace{1em}
\end{center}
\caption{Basic Planar 4R {\em Strange}  Linkage Types}
\label{Strange}
\end{table}

In the following we consider only normal mechanisms.

Here  we investigate first the rigid mechanisms.

\begin{thm} Suppose that $g=a+b+h$. Then 
the Minkowskian planar 4R has a rocker--crank type. 
\end{thm}
\proof

Since $T_1>0$, $T_2<0$, $T_3=0$, $T_4>0$ and $T_5>0$.

\qed
\begin{thm} Suppose that $a=g+b+h$. Then 
the Minkowskian planar 4R has a rocker--rocker type. 
\end{thm}
\proof

Since $T_1<0$, $T_2>0$, $T_3<0$, $T_4=0$ and $T_5>0$.

\qed

\begin{thm} Suppose that $h=g+b+a$. Then 
the Minkowskian planar 4R has a crank--rocker type. 
\end{thm}
\proof

Since $T_1<0$, $T_2<0$, $T_3<0$, $T_4>0$ and $T_5=0$.

\qed

\begin{thm} Suppose that $b=g+h+a$. Then 
the Minkowskian planar 4R has a crank--crank type. 
\end{thm}
\proof

Since $T_1>0$, $T_2>0$, $T_3<0$, $T_4>0$ and $T_5>0$.
\qed

We summarized the classification of rigid mechanisms in Table \ref{Rigid}.

\begin{table}[htp] 
\begin{center}
\begin{tabular}{||c|l|l||}
\hline
 & Linkage Type & Condition   \\
\hline
1 & Crank--Crank & $b=g+h+a$  \\
2 & Crank--Rocker & $h=g+b+a$   \\
3 & Rocker--crank & $g=a+b+h$ \\
4 & Rocker--Rocker & $a=g+b+h$ \\

\hline
\end{tabular}
\vspace{1em}
\end{center}
\caption{Basic Planar 4R {\em Rigid} Linkage Types}
\label{Rigid}
\end{table}

Now we consider normal, non--rigid reducible linkages.
\begin{thm} 
Suppose that $g+b=a+h$. Then the Minkowskian planar 4R has a crank--crank type.
\end{thm}
\proof 

Since $T_1=0$, hence $T_1T_2=0$. \qed

\begin{thm} 
Suppose that $a+b=g+h$. Then the input crank is crank.

\end{thm}
\proof 

Since $T_2=0$, hence $T_1T_2=0$. \qed

Finally we investigate non--rigid normal Minkowskian planar 4R mechanisms.
\begin{thm} 
For all non--rigid normal Minkowskian planar 4R we have $T_3<0$, $T_4>0$ and $T_5>0$.
\end{thm}
The link lengths $a,b,g$ and $h$ for a 4R chain define the two parameters 
$T_1$, $T_2$. Clearly our classification requires only the signs of these 
parameters. We assamble here an array for ($\mbox{sgn}T_1,\mbox{sgn}T_2$) (see Table \ref{Normal_Non_Rigid}).

%We can separate the linkage types into two general classes depending on the sign of
%the product $T_1T_2$. If $T_1T_2>0$, then the linkage is called {\em Grashof},
%otherwise it is called {\em nonGrashof}. 

\begin{table}[htp] 
\begin{center}
\begin{tabular}{||c|l|c|c||}
\hline
 & Linkage Type & $T_1$ & $T_2$  \\
\hline
1 & Crank--crank & + & +  \\
2 & Rocker--crank & + & --  \\
3 & Rocker--rocker & -- & +  \\
4 & Crank--rocker & -- & -- \\

\hline
\end{tabular}
\vspace{1em}
\end{center}
\caption{Basic Planar 4R {\em Normal Non--Rigid}  Linkage Types}
\label{Normal_Non_Rigid}
\end{table}

\begin{remark}
It is easy to see from the reversed triangle inequality (Lemma \ref{revers_spacelike})
that the following criterion is analogous to Grashof's criterion:
$$
l+s\geq p+q,
$$
where $s$ is  the length of the shortest link,
$l$ is the length of the longest link and $p$, $q$ are the lengths of the 
remaining two links.
  
%Hence there exists at least one link which can fully revolve with respect to
%the other three links if and only if 
%$$
%l+s\geq p+q.
%$$
\end{remark}

\subsection{Examples}

We animated the movement of the Minkowskian planar 4R in Matlab. In the following examples
 we show some pictures from our animation. 

1. Suppose that $a:=1$, $b:=1$, $g:=4$ and $h:=1$. Then $g>a+b+h$, hence this is a
{\em strange} Minkowskian planar 4R. We get
$$
\mbox{ch}(\theta_{min})= 1.625 ,
$$  
$$
\mbox{ch}(\theta_{max})= 2.125,
$$
$$
\mbox{ch}(\psi_{min})= -2.125,
$$
$$
\mbox{ch}(\psi_{max})= -1.625.
$$
$T_1=g+b-h-a=3$, $T_2=a-g+b-h=-3$, $T_3=g-a-b-h=1$, $T_4=g-a+b+h=5$ and $T_5=a-h+g+b=5$.
Since $T_1>0$, $T_2<0$, $T_3>0$, $T_4>0$ and $T_5>0$, 
thus this planar 4R has a superrocker--crank type.

The branching points occour at $\ch\theta =1.5$.

\begin{figure}[htp]
    \centering 
\includegraphics[height=180mm]{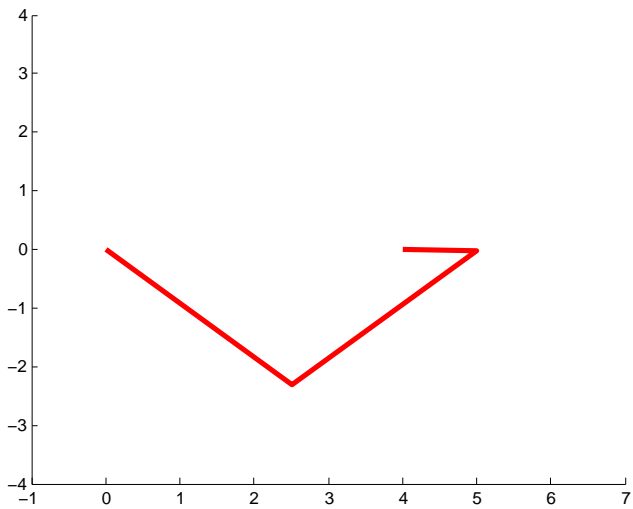}
\caption[Photo_st_min_pi_2]{The Minkowskian planar 4R with $a=1$, $b=1$, $h=1$, $g=4$ at $\theta=-\pi / 2$}
  \label{fig:Photo_st_min_pi_2}
\end{figure}

\begin{figure}[htp]
    \centering 
\includegraphics[height=180mm]{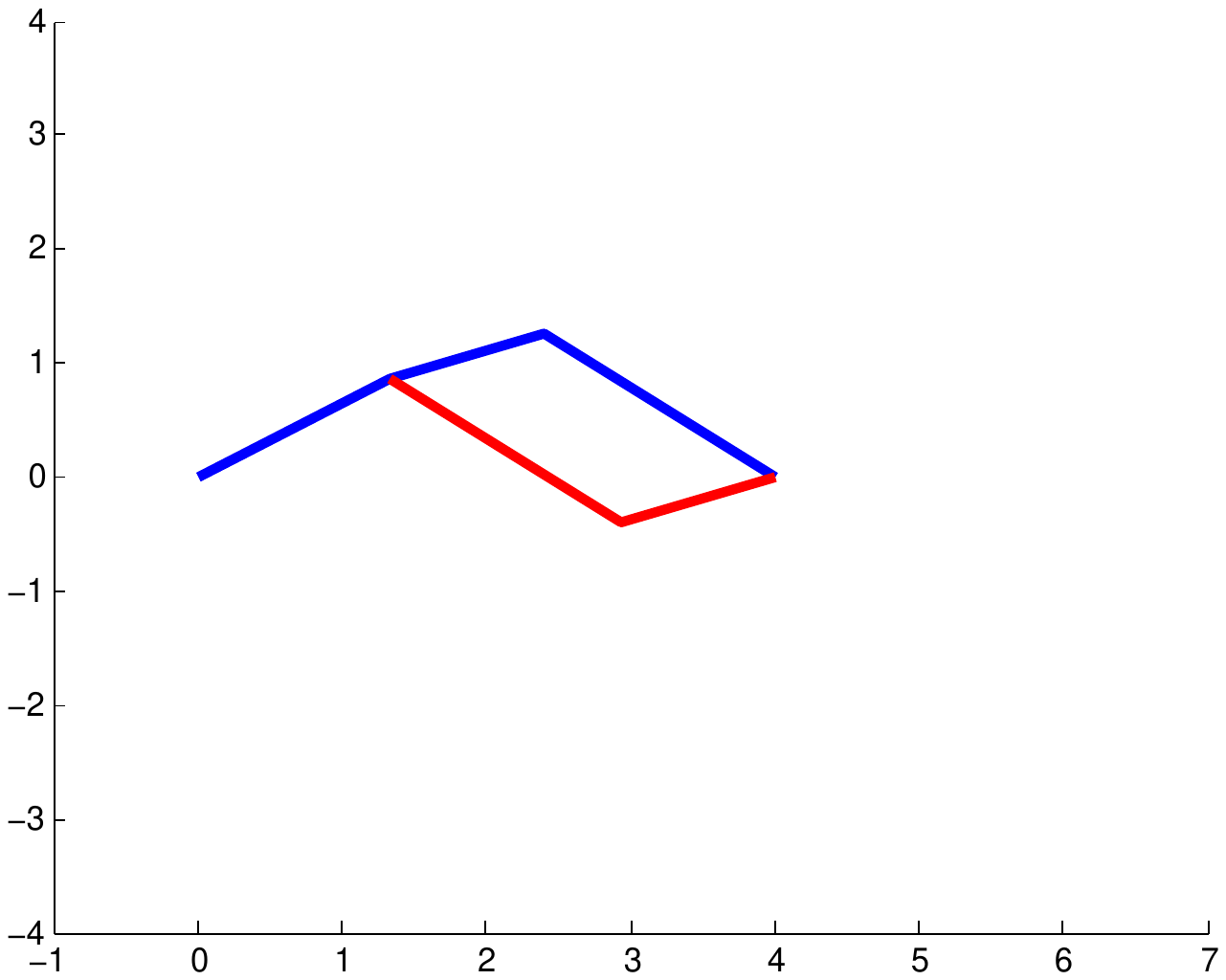}
\caption[Photo_st_pi_4]{The Minkowskian planar 4R with $a=1$, $b=1$, $h=1$, $g=4$ at $\theta=\pi / 4$}
  \label{fig:Photo_st_pi_4}
\end{figure}

\begin{figure}[htp]
    \centering 
\includegraphics[height=180mm]{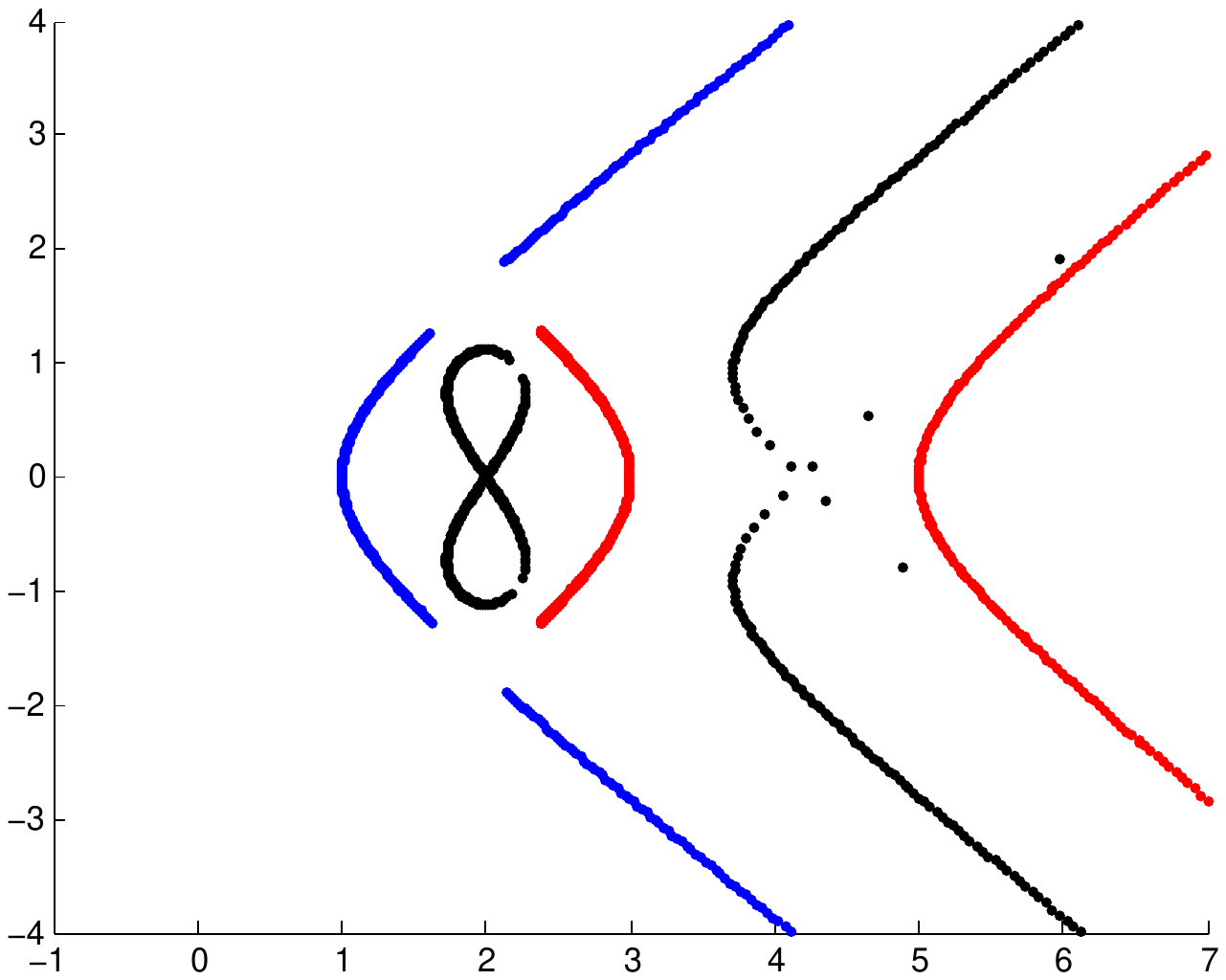}
\caption[coupler_st]{The coupler curve of the Minkowskian planar 4R with $a=1$, $b=1$, $h=1$, $g=4$}
  \label{fig:coupler_st}
\end{figure}

\begin{figure}[htp]
    \centering 
\includegraphics[height=180mm]{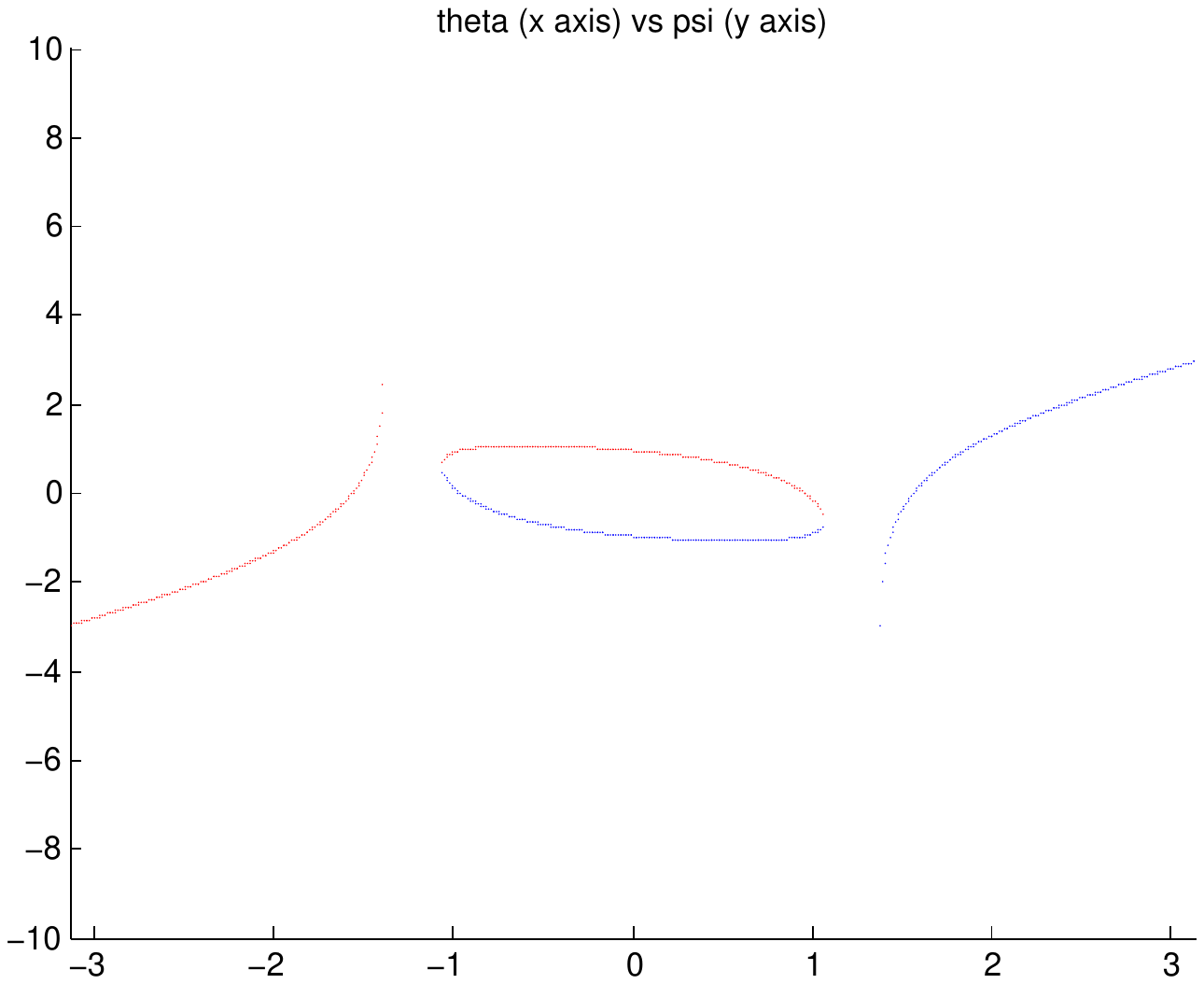}
\caption[trans_st]{The transmission function of the Minkowskian planar 4R with $a=1$, $b=1$, $h=1$, $g=4$}
  \label{fig:Trans_st}
\end{figure}

We show here two pictures from the animation (see Figures \ref{fig:Photo_st_min_pi_2}, 
 \ref{fig:Photo_st_pi_4}), the coupler curve (Figure \ref{fig:coupler_st}) and 
the transmission curve (Figure \ref{fig:Trans_st}). 
In Figure \ref{fig:coupler_st} the
blue curve shows the trajectory of the end of the input crank, 
the black curve shows the trajectory 
of the middle point of the coupler, while the red curve shows the 
trajectory of the end of the input crank


2. Suppose that $a:=1.2$, $b:=0.4$, $g:=0.4$ and $h:=0.4$. Then $a=g+b+h$, hence this is a
rigid Minkowskian planar 4R. We get
$$
\mbox{ch}(\theta_{min})= 1,
$$  
$$
\mbox{ch}(\theta_{max})\approx 1.6666,
$$
$$
\mbox{ch}(\psi_{min})= 1,
$$
$$
\mbox{ch}(\psi_{max})= 7.
$$
$T_1=g+b-h-a=-0.8$, $T_2=a-g+b-h=0.8$, $T_3=g-a-b-h=-1.6$, $T_4=g-a+b+h=0$ and $T_5=a-h+g+b=1.6$.
Since $T_1<0$, $T_2>0$, $T_3<0$, $T_4=0$ and $T_5>0$, 
thus this planar 4R has a rocker--rocker type.

Here $g=b$, but $h\neq a$, hence there are no branching points.

\begin{figure}[htp]
    \centering 
\includegraphics[height=180mm]{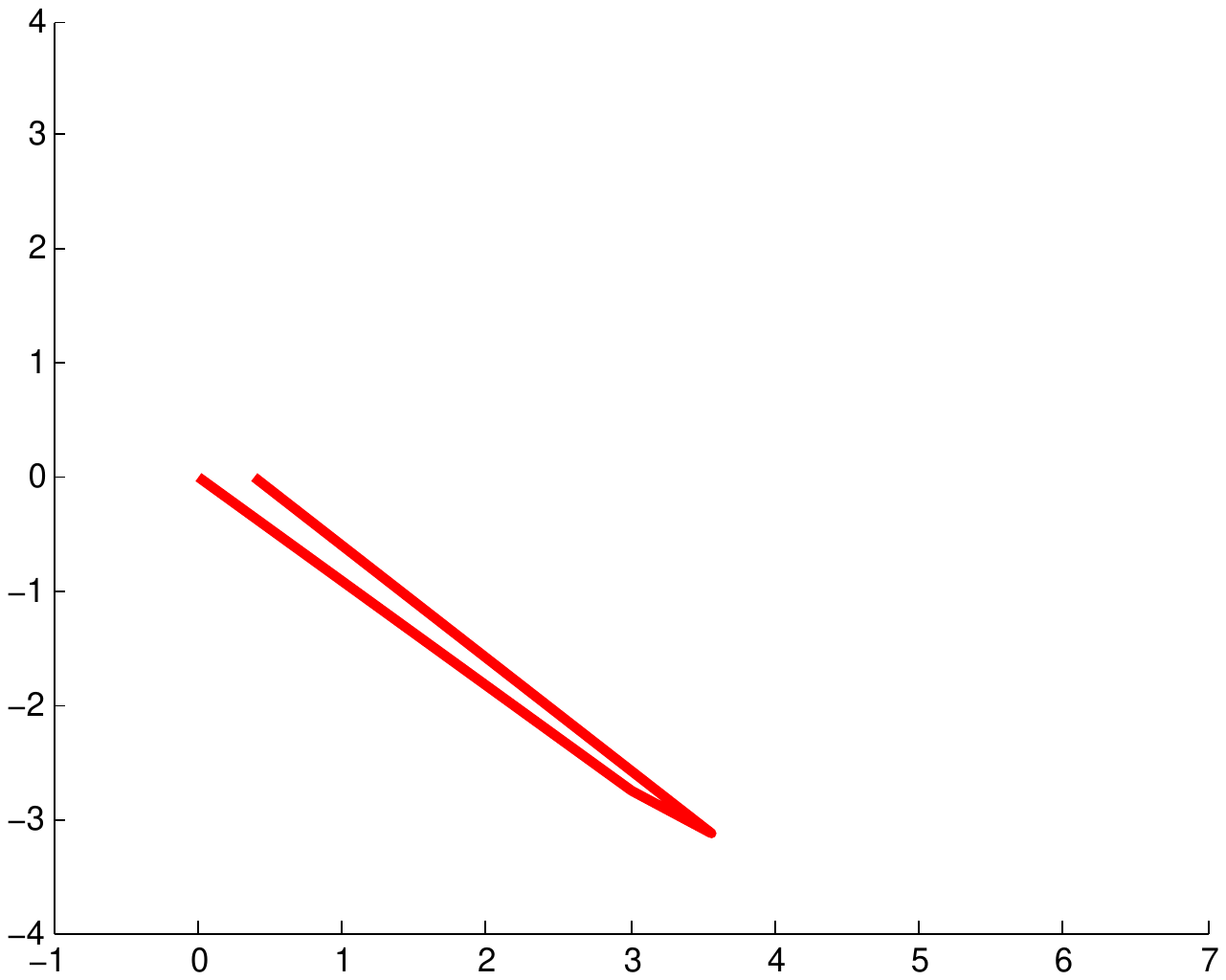}
\caption[Photo_rig_min_pi_2]{The Minkowskian planar 4R with $a=1.2$, $b=0.4$, $h=0.4$, $g=0.4$ at $\theta=-\pi / 2$}
  \label{fig:Photo_rig_min_pi_2}
\end{figure}

\begin{figure}[htp]
    \centering 
\includegraphics[height=180mm]{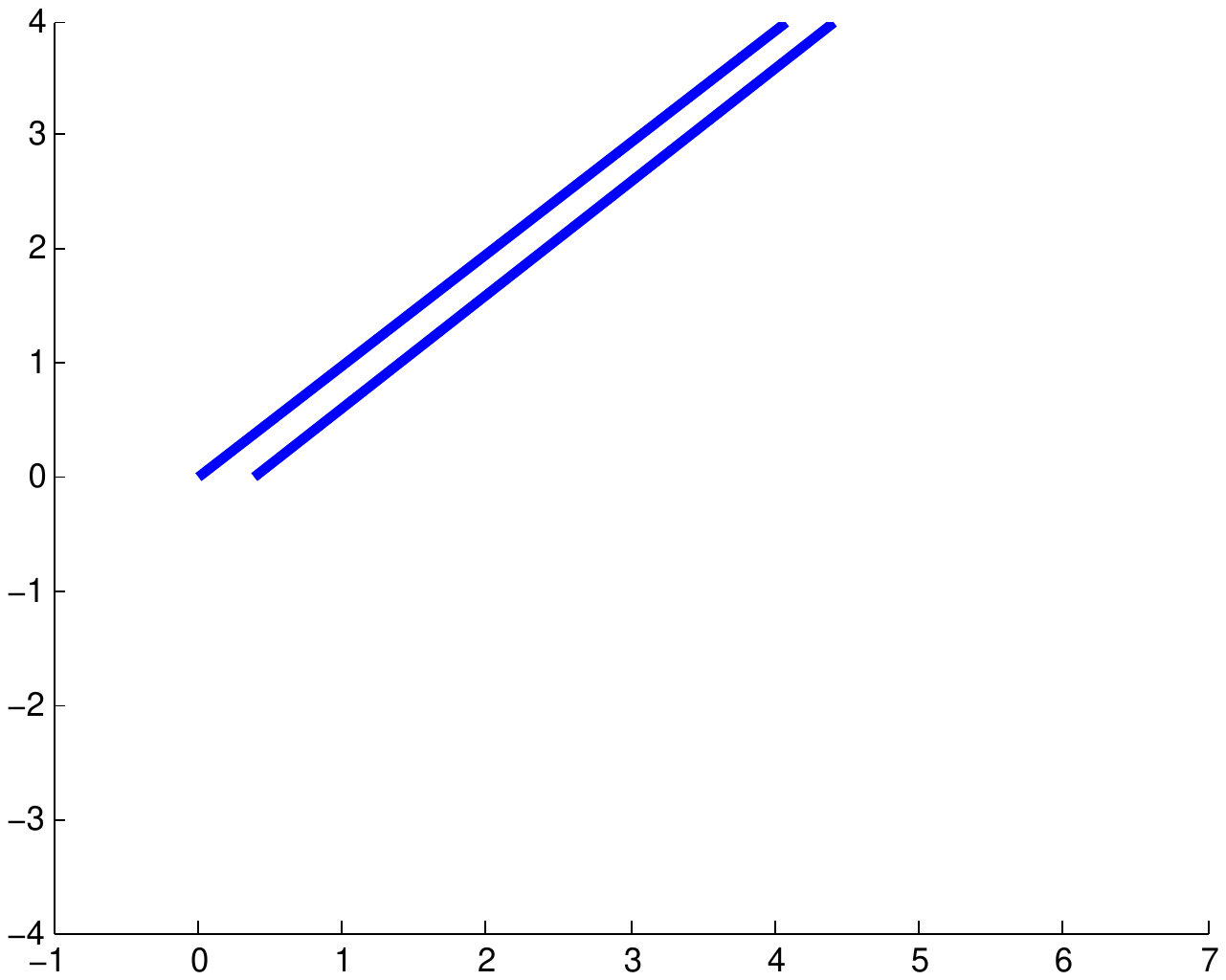}
\caption[Photo_rig_3_pi_4]{The Minkowskian planar 4R with $a=1.2$, $b=0.4$, $h=0.4$, $g=0.4$ at $\theta=3\pi / 4$}
  \label{fig:Photo_rig_3_pi_4}
\end{figure}

\begin{figure}[htp]
    \centering 
\includegraphics[height=180mm]{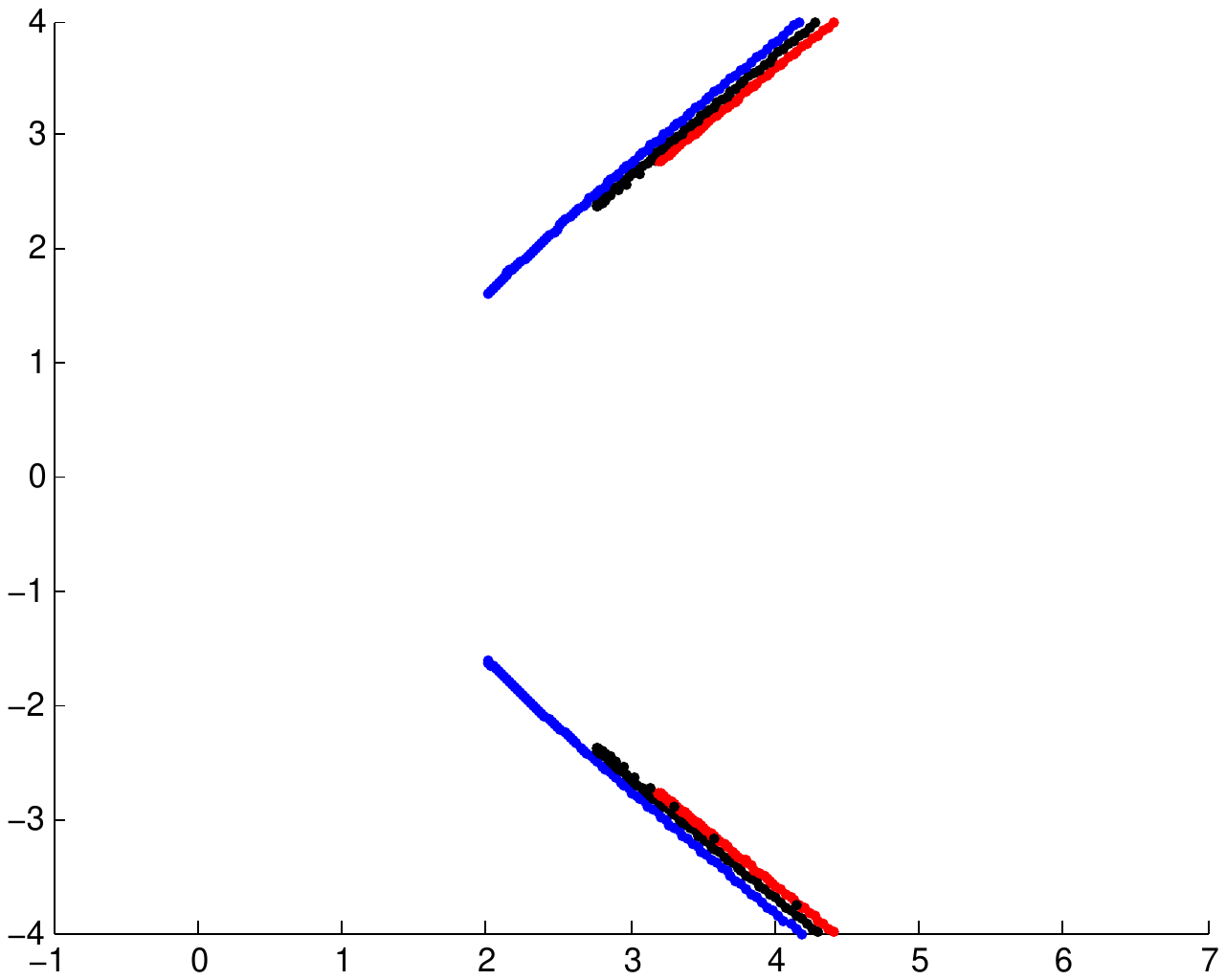}
\caption[coupler_rig]{The coupler curve of the Minkowskian planar 4R with $a=1.2$, $b=0.4$, $h=0.4$, $g=0.4$}
\label{fig:coupler_rig}
\end{figure}

\begin{figure}[htp]
    \centering 
\includegraphics[height=180mm]{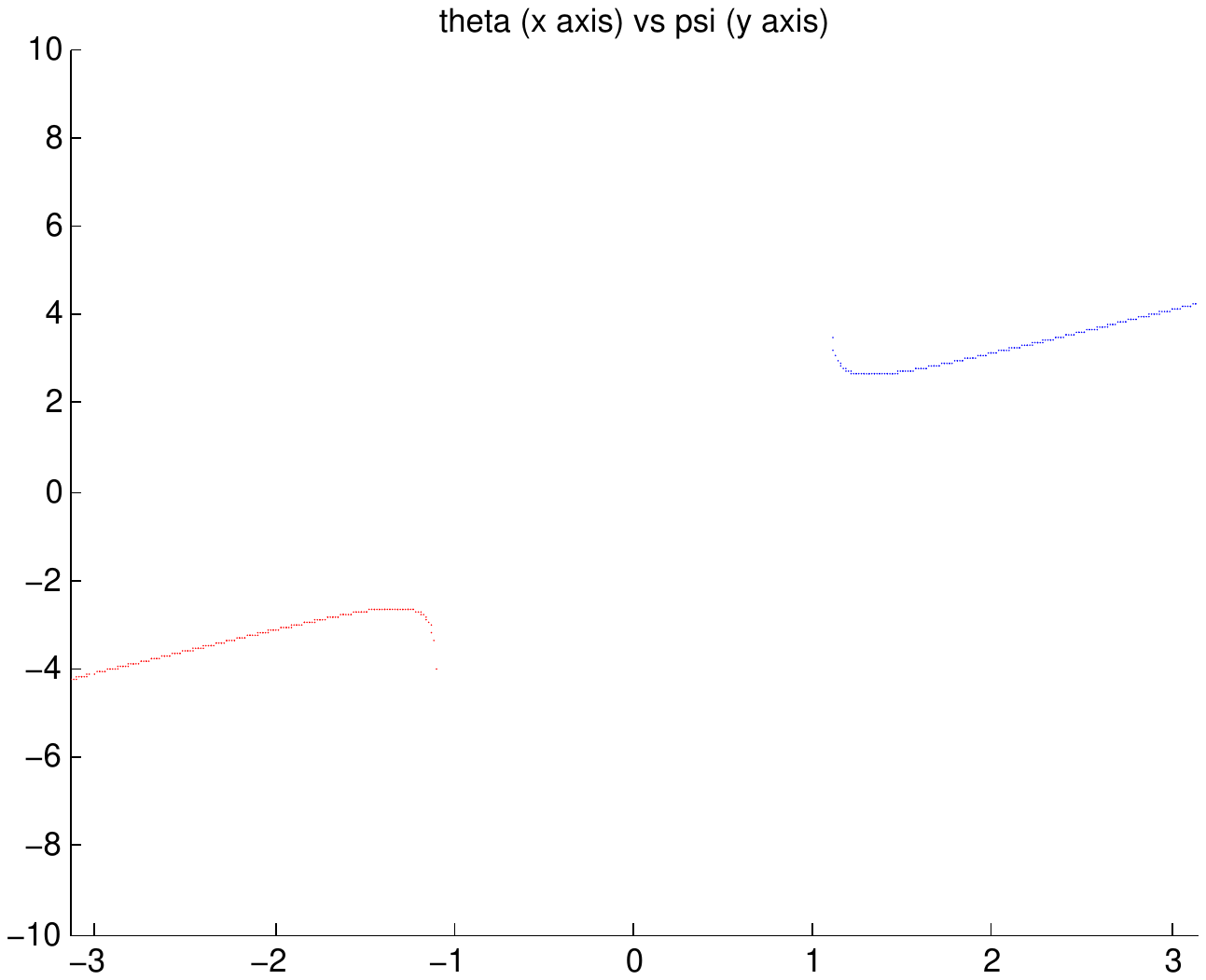}
\caption[trans_rig]{The transmission function of the Minkowskian planar 4R with $a=1.2$, $b=0.4$, $h=0.4$, $g=0.4$}
\label{fig:Trans_rig}
\end{figure}

We show here two pictures from the animation (Figures \ref{fig:Photo_rig_min_pi_2}, \ref{fig:Photo_rig_3_pi_4}), 
the coupler curve (Figure \ref{fig:coupler_rig}) and the transmission
curve (Figure \ref{fig:Trans_rig}).
In Figure \ref{fig:coupler_rig} the
blue curve shows the trajectory of the end of the input crank, 
the black curve shows the trajectory 
of the middle point of the coupler, while the red curve shows the 
trajectory of the end of the input crank


3. Suppose that $a:=0.5$, $b:=1$, $g:=2$ and $h:=2.5$. Then $g+b=a+h$, hence this is a
normal, non--rigid reducible Minkowskian planar 4R. We get
$$
\mbox{ch}(\theta_{min})= -4,
$$  
$$
\mbox{ch}(\theta_{max})= 1,
$$
$$
\mbox{ch}(\psi_{min})= -0.25,
$$
$$
\mbox{ch}(\psi_{max})= 1.
$$
Hence $T_1=g+b-h-a=0$, $T_2=a-g+b-h=-3$, $T_3=g-a-b-h=-2$, $T_4=g-a+b+h=5$ and $T_5=a-h+g+b=1$.
Since $T_1=0$, $T_2<0$, $T_3<0$, $T_4>0$ and $T_5>0$, 
thus this planar 4R has a crank--crank type.

The branching points occour at $\ch\theta =-5$, hence there are no branching points..

\begin{figure}[htp] 
    \centering 
\includegraphics[height=180mm]{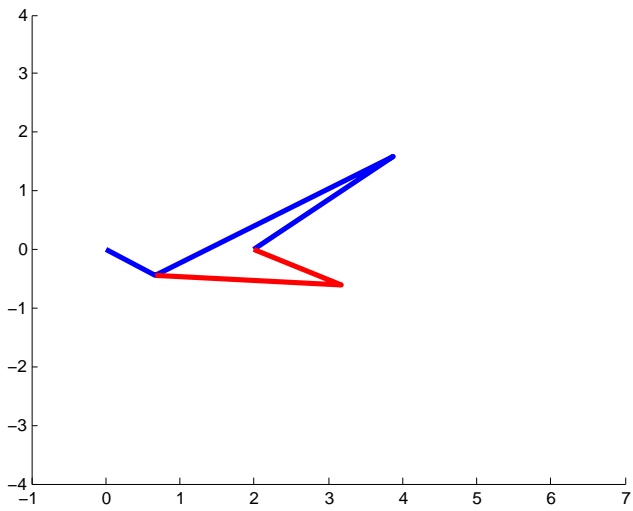}
\caption[Photo_red_min_pi_4]{The Minkowskian planar 4R with $a:=0.5$, $b:=1$, $g:=2$, $h:=2.5$ at $\theta=-\pi / 4$}
 \label{fig:Photo_red_min_pi_4}
\end{figure}

\begin{figure}[htp]
    \centering 
\includegraphics[height=180mm]{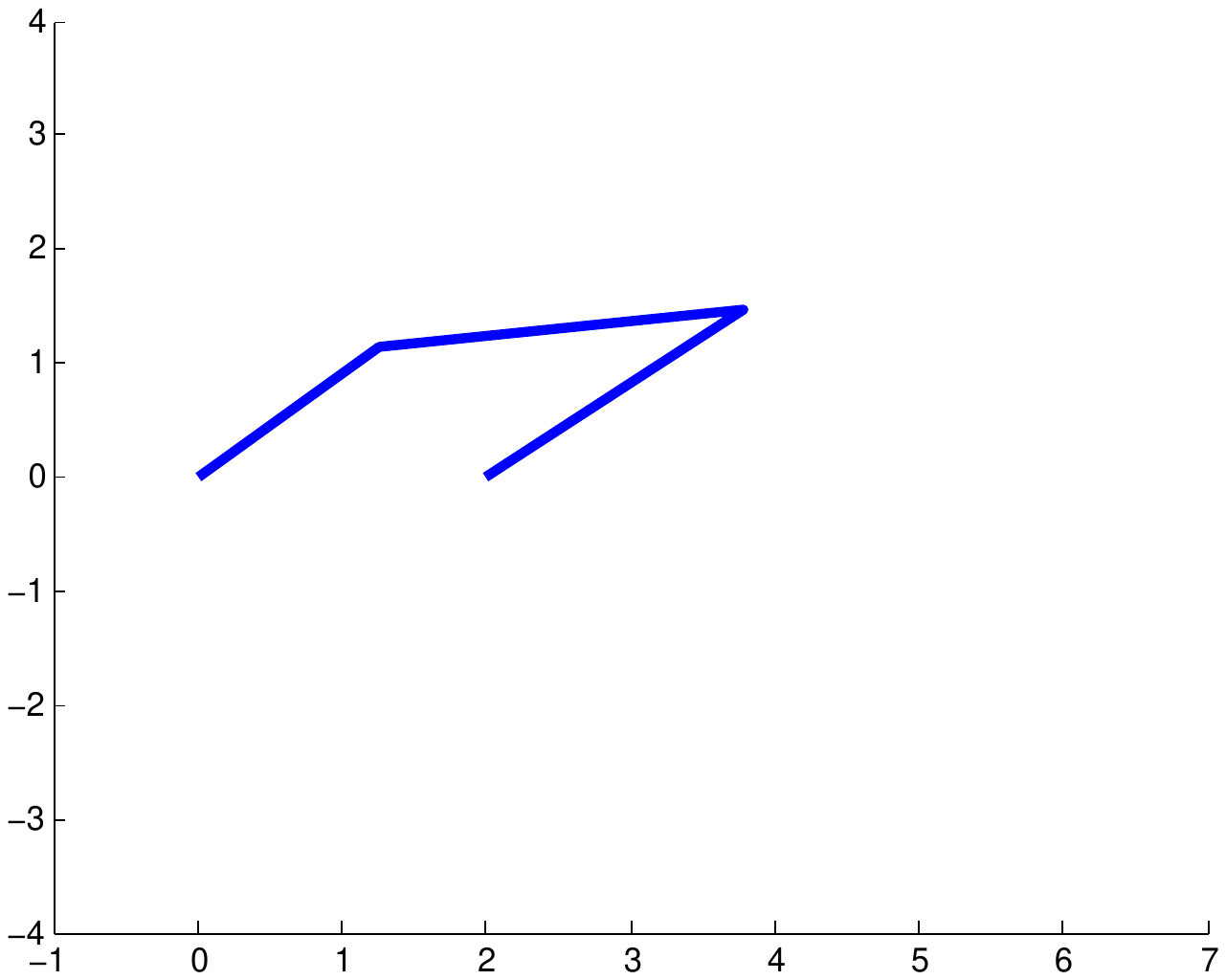}
\caption[Photo_red_pi_2]{The Minkowskian planar 4R with $a:=0.5$, $b:=1$, $g:=2$, $h:=2.5$ at $\theta=\pi / 2$}
\label{fig:Photo_red_pi_2}
\end{figure}

\begin{figure}[htp] 
    \centering 
\includegraphics[height=180mm]{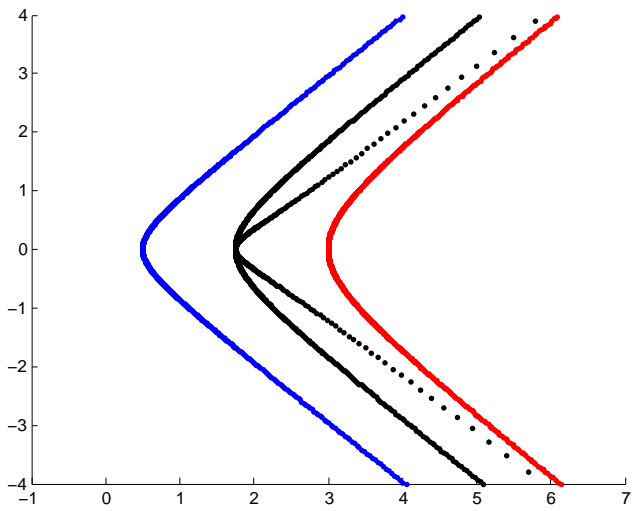}
\caption[coupler_red]{The coupler curve of the Minkowskian planar 4R with $a:=0.5$, $b:=1$, $g:=2$, $h:=2.5$}
\label{fig:coupler_red}
 \end{figure}

\begin{figure}[htp]
    \centering 
\includegraphics[height=180mm]{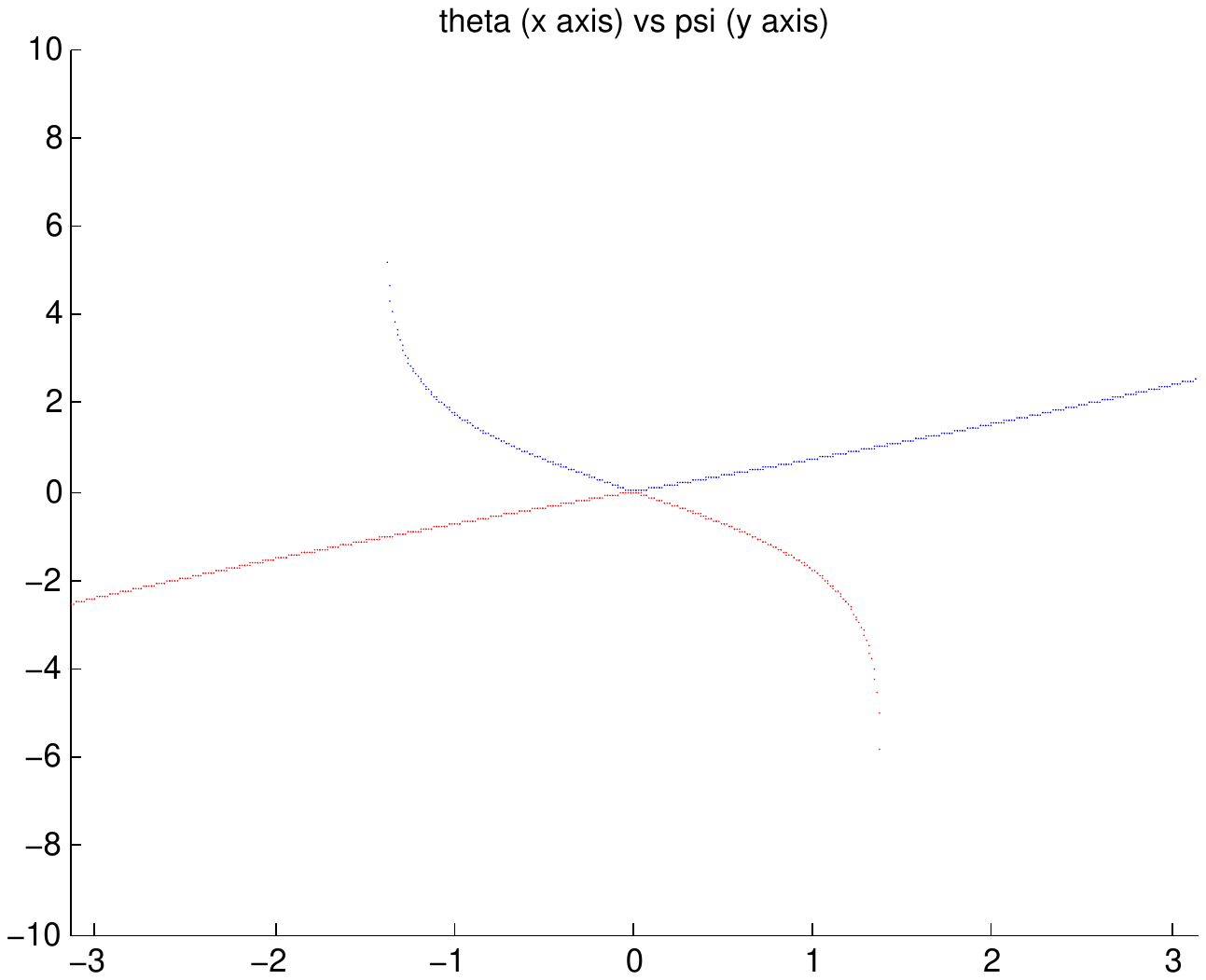}
\caption[trans_red]{The transmission function of the Minkowskian planar 4R with $a:=0.5$, $b:=1$, $g:=2$, $h:=2.5$}
\label{fig:Trans_red}
 \end{figure}

We show here two pictures from the animation (Figures \ref{fig:Photo_red_min_pi_4}, \ref{fig:Photo_red_pi_2}), 
the coupler curve (Figure \ref{fig:coupler_red}) and 
the transmission curve (Figure \ref{fig:Trans_red}).
In Figure \ref{fig:coupler_red} the
blue curve shows the trajectory of the end of the input crank, 
the black curve shows the trajectory 
of the middle point of the coupler, while the red curve shows the 
trajectory of the end of the input crank

4. Suppose that $a:=0.6$, $b:=1$, $g:=0.7$ and $h:=0.5$. Then this is a
normal, non--rigid irreducible Minkowskian planar 4R. We get
$$
\mbox{ch}(\theta_{min})\approx -1.6666,
$$  
$$
\mbox{ch}(\theta_{max})\approx 0.71428,
$$
$$
\mbox{ch}(\psi_{min})\approx -1.05714,
$$
$$
\mbox{ch}(\psi_{max})= -0.2.
$$
Hence $T_1=g+b-h-a=0.6$, $T_2=a-g+b-h=0.4$, $T_3=g-a-b-h=-1.4$, $T_4=g-a+b+h=1.6$ and $T_5=a-h+g+b=1.8$.
Since $T_1>0$, $T_2>0$, $T_3<0$, $T_4>0$ and $T_5>0$, 
thus this planar 4R has a crank--crank type.

The branching points occour at $\ch\theta \approx -0.5555$, hence there are no branching points..

\begin{figure}[htp] 
    \centering 
\includegraphics[height=180mm]{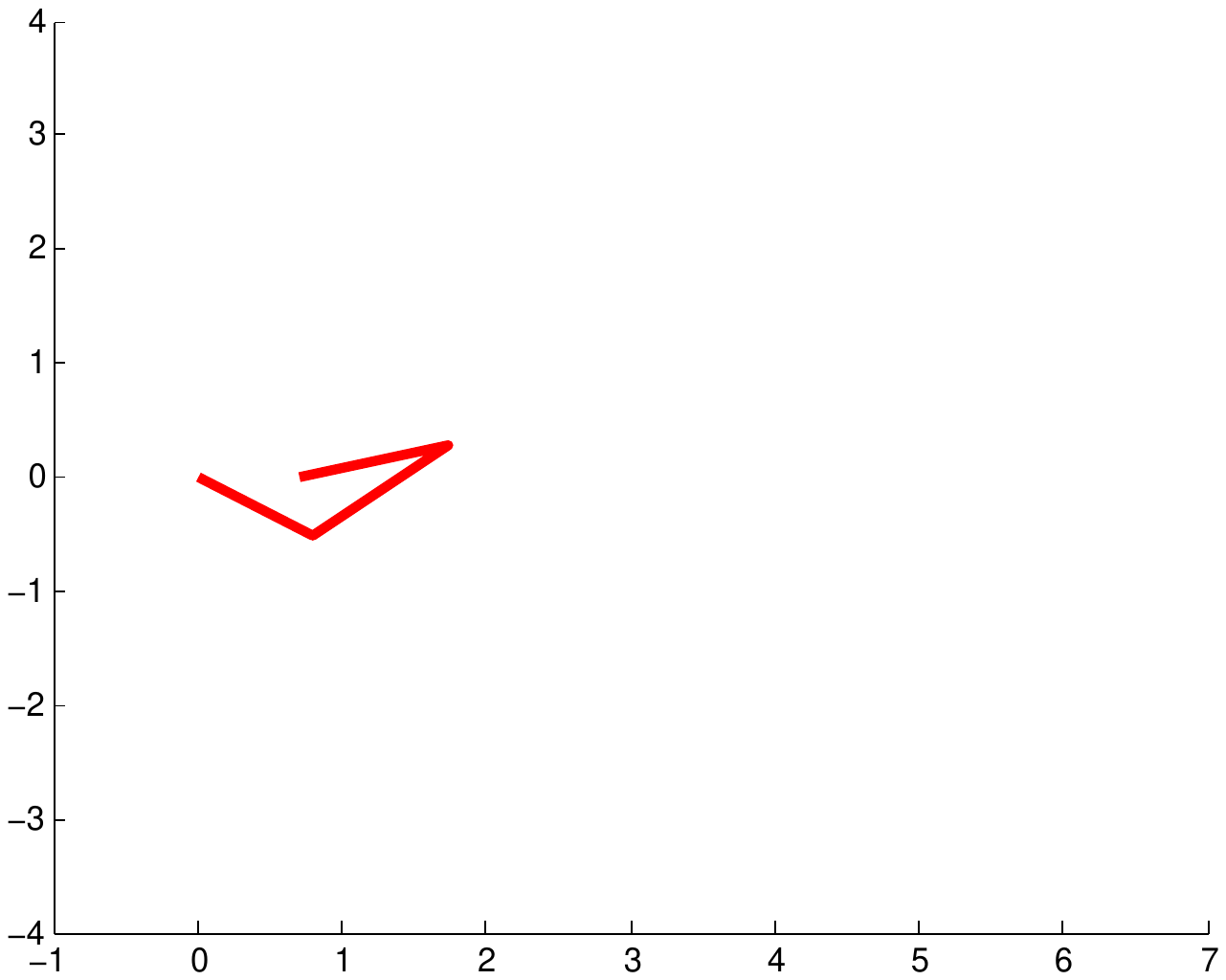}
\caption[Photo_irred_min_pi_4]{The Minkowskian planar 4R with $a:=0.6$, $b:=1$, $g:=0.7$, $h:=0.5$ at $\theta=-\pi / 4$}
  \label{fig:Photo_irred_min_pi_4}
\end{figure}

\begin{figure}[htp] 
    \centering 
\includegraphics[height=180mm]{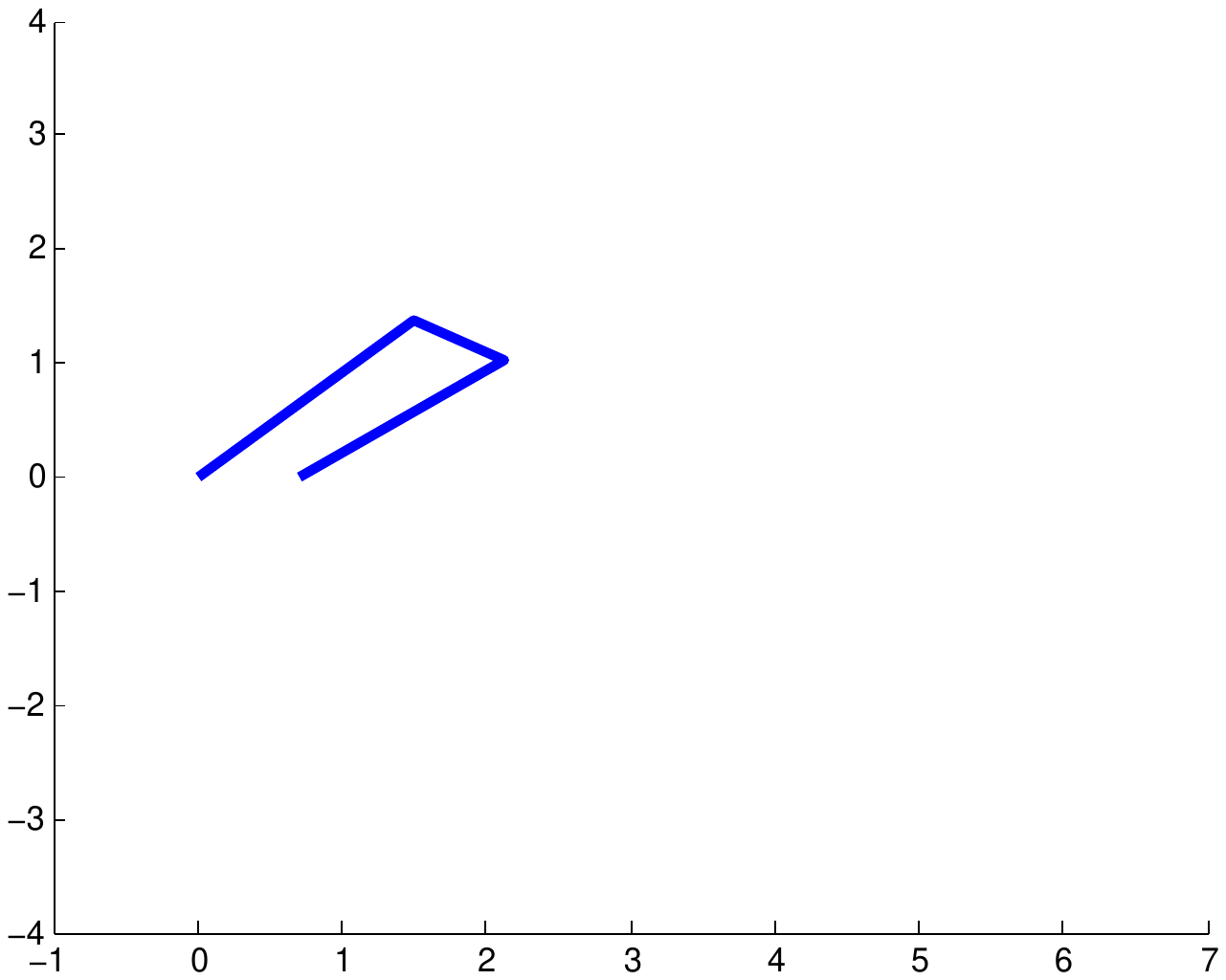}
\caption[Photo_irred_pi_2]{The Minkowskian planar 4R with $a:=0.6$, $b:=1$, $g:=0.7$, $h:=0.5$ at $\theta=\pi / 2$}
\label{fig:Photo_irred_pi_2}
\end{figure}

\begin{figure}[htp] 
    \centering 
\includegraphics[height=180mm]{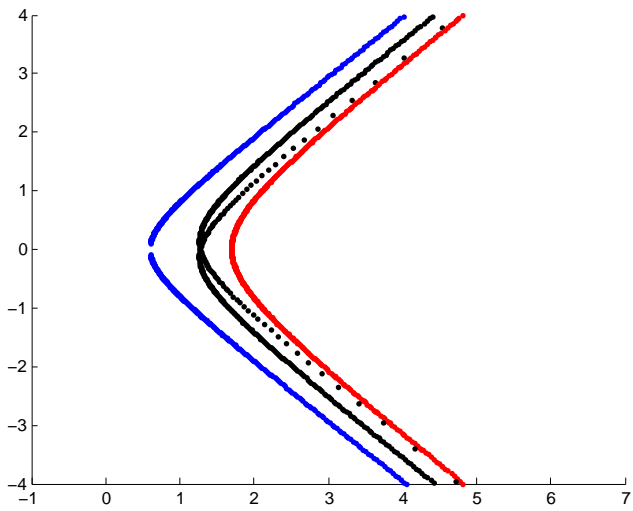}
\caption[coupler_irred]{The coupler curve of the Minkowskian planar 4R with $a:=0.6$, $b:=1$, $g:=0.7$, $h:=0.5$}
 \label{fig:coupler_irred}
\end{figure}

\begin{figure}[htp] 
    \centering 
\includegraphics[height=180mm]{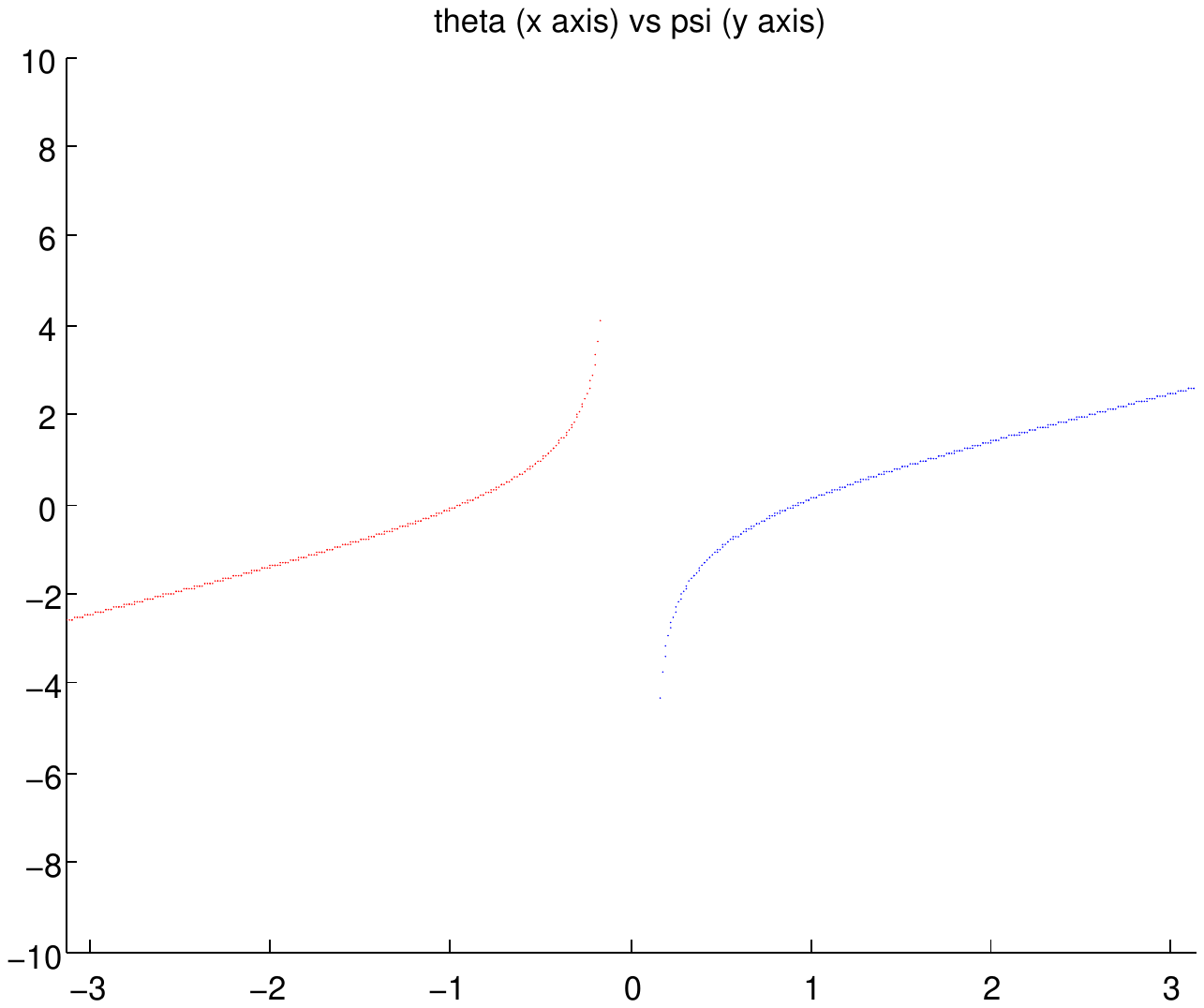}
\caption[trans_red]{The transmission function of the Minkowskian planar 4R with $a:=0.6$, $b:=1$, $g:=0.7$, $h:=0.5$}
\label{fig:Trans_irred}
\end{figure}

We show here two pictures from the animation (Figures  \ref{fig:Photo_irred_min_pi_4},\ref{fig:Photo_irred_pi_2}), 
the coupler curve (Figure \ref{fig:coupler_irred}) and 
the transmission curve (Figure  \ref{fig:Trans_irred}).
In Figure \ref{fig:coupler_irred} the
blue curve shows the trajectory of the end of the input crank, 
the black curve shows the trajectory 
of the middle point of the coupler, while the red curve shows the 
trajectory of the end of the input crank

\section{Conclusion}

In this article we characterized and classified completely 
the planar 4R
closed chain working on the Minkowskian plane. We derived formulas for the 
output  crank angle, the coupler angle and the transmission 
angle. We found four basic types in the classification: the crank-crank, 
crank--rocker, rocker--crank and rocker--rocker Minkowskian planar 4R mechanisms and 
we described the Minkowskian Grashof condition.

{\bf Acknowledgements.} The first author is grateful to Josef Schicho and Madalina Hodorog for their useful
comments.

\end{document}